\documentclass[reqno]{amsart}

\usepackage[foot]{amsaddr}
\usepackage{comment}
\usepackage[utf8]{inputenc}
\usepackage{graphicx}
\usepackage{verbatim}
\usepackage{amsmath}
\usepackage{mathtools}
\usepackage{mathabx}
\usepackage{siunitx}
\usepackage[dvipsnames, table]{xcolor}
\usepackage[top=1in, bottom=1.25in, left=1.25in, right=1.25in]{geometry}
\usepackage{hyperref}
\usepackage{algorithm}
\usepackage{algpseudocode}

\usepackage{ulem}
\usepackage{enumitem}

\usepackage{booktabs}
\usepackage{multirow}
\usepackage{multicol}
\usepackage{pifont}

\usepackage{subfig}

\def\PP{{{\rm l}\kern - .15em {\rm P} }}
\def\PN2{{\PP_{N}-\PP_{N-2}}}

\newcommand{\bu}{\boldsymbol{u}}

\newcommand{\bw}{\boldsymbol{w}}

\newcommand{\bx}{\boldsymbol{x}}

\newcommand{\deleted}[1]{{}}

\theoremstyle{remark}
\newtheorem{remark}{Remark}[section]
\title[RL filters with and without reference data training]{Reinforcement Learning-Based Filters for Convection-Dominated Flows: Reference-Free and Reference-Guided Training}

\author{Anna Ivagnes$^{1, *}$, Maria Strazzullo$^{2, *}$ and Gianluigi Rozza$^{1, *}$}
\address{$^1$ mathLab, Mathematics Area, SISSA, via Bonomea 265, I-34136 Trieste, Italy}
\address{$^2$ Politecnico di Torino, Department of Mathematical Sciences ``Giuseppe Luigi Lagrange'', Corso Duca degli Abruzzi, 24, 10129, Torino, Italy}
\address{$^*$ INdAM-GNCS group member}

\begin{document}
\begin{abstract}
We propose a reinforcement learning (RL) framework for the dynamic selection of the filter parameter in Evolve–Filter (EF) regularization strategies for incompressible turbulent flows. Instead of prescribing the filter radius heuristically, the RL agent learns to adaptively control the filtering intensity in time, balancing numerical stability and physical accuracy.

The methodology is assessed on two benchmark problems with fundamentally different dynamics: flow past a cylinder and decaying homogeneous turbulence. Both reference-guided and reference-free reward formulations are investigated. In the reference-guided setting, the agent is trained using direct numerical simulation (DNS) data over a limited time window and then evaluated in extrapolation. In the reference-free setting, the reward relies exclusively on physics-based quantities, without access to reference solutions, i.e., completely eliminating the computational costs related to DNS simulations.

The results show that the proposed RL-EF strategies prevent numerical blow-up while avoiding the excessive dissipation typical of standard EF approaches based on a fixed Kolmogorov length scale. The learned policies accurately reproduce the relevant flow dynamics across scales, preserving the correct balance between large-scale structures and small-scale dissipation. Notably, the reference-free reward achieves performance comparable to the reference-guided approach, demonstrating that stable and spectrally consistent filtering strategies can be learned even without DNS data, drastically reducing the computational costs of the training phase.

The proposed framework provides a robust and flexible alternative to manually tuned regularization parameters, enabling adaptive, physically consistent control of filtering in turbulent flow simulations.
\end{abstract}

\maketitle

\section{Introduction}
Convection-dominated flows pose significant challenges in under- or marginally-resolved regimes simulations, i.e., when using computational grids that do not scale down to the Kolmogorov scale. In such settings, classical discretizations of the Navier–Stokes Equations, such as finite-element or finite-volume strategies, may produce spurious oscillations, loss of stability, and nonphysical energy accumulation~\cite{berselli2006mathematics,layton2012approximate}. To address these issues, stabilization mechanisms are usually incorporated, including large eddy simulation (LES)~\cite{choi2025performance, li2025optimizing}, Reynolds-averaged Navier–Stokes (RANS) modeling~\cite{berselli2006mathematics,layton2012approximate,rebollo2014mathematical,guo2025effects, ahmadi2024predicting, martinez2025rans, zou2025study}, and regularization-based approaches. Regularization (see \cite{layton2012approximate} for an overview of these methods) is a class of numerical schemes that exploit spatial filtering to increase the accuracy of numerical flow simulations, allowing the employment of less refined mesh-based discretization if compared to the direct numerical simulation (DNS), maintaining accuracy while reducing the computational costs.
Among the regularization approaches, the Evolve–Filter (EF) strategy has emerged as an effective and computationally inexpensive framework. The key idea is to repeat at each time step of the simulation: (i) an \textit{evolve} step, e.g., advance the velocity with standard Navier--Stokes discretization; (ii) a spatial \textit{filtering} stage, which smooths out the evolved velocity and alleviates the spurious oscillations.
A widely used class of spatial filters is the differential filter (DF). The success of DFs 
is due to their mathematical robustness \cite{berselli2006mathematics}.

EF-type methods can be easily integrated into existing solvers, which explains their widespread adoption across discretization paradigms~\cite{boyd2001chebyshev, canuto1988some,ervin2012numerical,fischer2008nek5000,fischer2001filter,germano1986differential, germano1986differential2,girfoglio2019finite, girfoglio2021fluid,xie2018evolve,mullen1999filtering,NekROM,olshanskii2013connection,pasquetti2002comments, rezaian2023predictive,takhirov2018computationally,tsai2025time,strazzullo2023new,strazzullo2022consistency} and a wide range of applications~\cite{girfoglio2019finite, girfoglio2023novel, besabe2025linear, xuthesis, quaini2024bridging, xu2020backflow,tsai2023parametric,layton2012approximate}.
However, a central issue in EF formulations is the selection of the filter radius \cite{bertagna2016deconvolution,girfoglio2019finite,strazzullo2022consistency, strazzullo2024variational,tsai2025time, berselli2006mathematics,ZoccolanStrazzulloRozza20242,ZoccolanStrazzulloRozza2024}. The method is highly sensitive to this parameter: insufficient filtering may lead to instability, whereas excessive filtering results in over-diffusive solutions and loss of physically relevant structures. Most existing approaches prescribe constant parameter values both in space and time based on heuristic arguments, such as the mesh size or the Kolmogorov length scale. Recent works have considered space-dependent parameters via nonlinear indicator functions in filters~\cite{girfoglio2019finite}, and time-dependent parameters via optimization with respect to reference DNS data. However, the first case poses the problem of a proper selection of the indicator function, while the second case heavily relies on data (usually related to very high computational costs) at each time instance and does not permit generalization beyond the calibration window.
In this work, we reformulate the parameter optimization as a sequential decision-making problem.
Instead of optimizing the filter parameter independently at each time step, we learn a time-dependent control policy through reinforcement learning (RL) strategies.

RL represents a class of machine learning methods in which an agent interacts with an environment and learns a policy by maximizing a cumulative reward signal~\cite{wiering2012reinforcement, kaelbling1996reinforcement}. In contrast to supervised learning, RL does not require labeled input–output pairs; rather, it discovers optimal actions through sequential trial-and-error interactions. This paradigm is particularly well suited for control and decision-making problems where actions influence future states, making it a natural framework for time-dependent parameter selection.
In computational fluid dynamics, RL has recently been employed in flow control, turbulence modeling, subgrid-scale closure modeling, and adaptive numerical stabilization \cite{kurz2022deep, kurz2023deep, kurz2025harnessing, rabault2019artificial, rabault2020deep, tang2020robust, font2025deep, fan2020reinforcement, vinuesa2022flow, ren2021applying, paris2021robust}. These applications demonstrate its capability to learn nontrivial control strategies directly from high-dimensional flow states, even in complex nonlinear regimes. In this work, we leverage these strengths to dynamically regulate the filtering intensity within the EF framework, casting parameter selection as an adaptive control problem driven by the evolving flow dynamics.

While in previous works the agent training stage in reinforcement learning mostly relies on reference data, we here consider both DNS reference-driven and DNS reference-free formulations. In the following, we refer to the reference-free RL-based filter formulation as data-free (DF-RL-EF), meaning that no reference DNS data are used in the reward construction, in contrast with the reference data-driven setting (DD-RL-EF).  
Specifically, in the DD-RL-EF setting, the reward incorporates reference data over a limited training interval, and the learned policy is subsequently evaluated in extrapolation. In the DF-RL-EF, the reward is constructed exclusively from physically meaningful quantities, without access to reference solutions. This formulation allows for flexible reward design, including the incorporation of global invariants or structure-preserving constraints.
The resulting framework replaces static parameter calibration with adaptive, state-dependent control. Numerical experiments on representative convection-dominated flows demonstrate that the learned policies achieve a favorable balance between numerical stability and physical accuracy, while exhibiting robustness across distinct dynamical regimes.
The main novelties of this contribution rely on:
\begin{itemize}

\item  the introduction of a novel RL-selected filter radius selection for the EF strategy, guaranteeing an adaptive selection of the filtered procedure non-trivially in time, preserving large-scale and maintaining the main features of the flow simulation if compared to the reference data;

\item
Most importantly, the investigation of an RL-EF algorithm that does not rely on the DNS data, based on physics quantities such as equation residual, gradient-based or structure-preserving indicators, capable of providing a stable and reliable filter action with a much reduced computational cost.
\end{itemize}
The paper is organized as follows. Section \ref{sec:ef} reviews the EF formulation and discusses parameter sensitivity. Section \ref{sec:rl-ef-all} introduces the RL theory, the RL-EF methodology, and the proposed reward formulations. Section \ref{sec:results} presents the numerical investigations on two-dimensional fluid dynamics benchmarks: (i) the flow past a cylinder at $Re=\num{1000}$, and (ii) the decaying homogeneous turbulence at $Re=\num{40000}$. Conclusions and future research directions are outlined in Section \ref{sec:conclusions}.

\section{Problem Formulation and Evolve-Filter Strategy}
\label{sec:ef}
In this section we describe the backbone of the methodology we investigate.
Specifically, we briefly introduce the discrete Navier--Stokes Equations (NSE) and the Evolve and Filter (EF) strategy in Section 
\ref{subsec:NSE}. The interested reader may refer to \cite{wells2017evolve, strazzullo2022consistency, girfoglio2021fluid,bertagna2016deconvolution} for more details on EF.

\subsection{The EF algorithm for NSE}
\label{subsec:NSE}
In the proposed numerical investigation, we focus on the \textit{incompressible NSE}. At the continuous level, for a physical domain $\Omega \subset \mathbb R^{d}$, with $d=2,3$, we model the motion of an incompressible fluid of velocity $\bu = \bu(\bx, t) \in \mathbb U$ and pressure $p = p(\bx, t) \in \mathbb Q$ by means of the following equations:
\begin{equation}
\label{eq:NSE}
\begin{cases}
 \displaystyle \frac{\partial \bu}{\partial t} + (\bu \cdot \nabla) \bu - \nu \Delta \bu + \nabla p = 0 & \text{in }  \Omega \times (t_0, T), \\
\nabla \cdot \bu = 0  & \text{in }  \Omega \times (t_0, T), \\
\bu = \bu_D & \text{on } \Gamma_D \times (t_0, T), \\
\displaystyle -p \boldsymbol n + \nu \frac{\partial \boldsymbol \bu}{\partial \boldsymbol n} = 0  & \text{on }  \Gamma_N \times (t_0, T), \\
\end{cases}
\end{equation}
for a given initial condition $\bu = \bu_0$ in $\Omega \times \{ t_0 \}$, with $\Gamma_D \cup \Gamma_N = \partial \Omega$, $\Gamma_D \cap \Gamma_N = \emptyset$, and where $\nu$ is the kinematic viscosity, with $\mathbb U$ and $\mathbb Q$ are suitable Hilbert function spaces. On $\Gamma_D$, the given function $\bu_D$ is prescribed as a Dirichlet boundary condition, while on $\Gamma_N$ ``free flow" boundaries are applied. 
The flow regime is determined by the adimensional \textit{Reynolds number}, defined as 
\begin{equation}
\label{eq:Re}
Re = \frac{UL}{\nu},
\end{equation}
where $U$ and $L$ are the characteristic velocity and length of the problem at hand, respectively. 
For a large $Re$, the inertial forces dominate the viscous forces and the flow regime is referred to as \textit{convection-dominated}.
In the proposed numerical tests, we employ the $\mathbb P^2$-$\mathbb P^1$ Taylor-Hood \textit{Finite Element Method (FEM)} and a backward differentiation formula of order 1 (BDF1) for space and time discretization, respectively. When working with under-resolved or marginally-resolved numerical simulations, standard spatial discretization techniques might yield spurious numerical oscillations.
To tackle this issue, we employ EF-based algorithms. For a given time step $\Delta t$, we consider the time instances $t_n =  t_0 + n\Delta t$ for $n = 0, \dots, N$, with final time $T = N\Delta t$. At time $t_n$, $\bu{_n} \in \mathbb U^{N_h^{\bu}}$ and $p{_n} \in \mathbb Q^{N_h^p}$ represent the semi-discrete FEM velocity and pressure solutions, respectively, where ${N_h^{\bu}}$ and ${N_h^p}$ are the degrees of freedom related to the FEM spaces approximation $\mathbb U^{N_h^{\bu}}$ and $\mathbb Q^{N_h^{p}}$.
The EF algorithm at time $t_{n+1}$ for a pre-determined filter radius $\delta$ reads:
\begin{eqnarray*}      
&	\text{\bf (I)}& \text{\emph{Evolve}:} \quad 
\begin{cases}
\begin{split}
     &\frac{\bw_{n + 1} - \bu_n}{\Delta t} + (\bw_{n+1} \cdot \nabla) \bw_{n+1} -\nu \Delta \bw_{n+1} + \nabla p_{n+1} = 0 
\end{split} & \text{in } \Omega \times \{t_{n+1}\}, \vspace{2mm}\\
\nabla \cdot \bw_{n+1} = 0 & \text{in } \Omega \times \{t_{n+1}\}, \vspace{1mm}\\
\bw_{n+1} = \bu_D & \text{on } \Gamma_D \times \{t_{n+1}\}, \vspace{2mm}\\
\displaystyle -p_{n+1} \cdot \boldsymbol n + \frac{\partial \boldsymbol \bw_{n+1}}{\partial \boldsymbol n} = 0  & \text{on } \Gamma_N \times \{t_{n+1}\}. \\
\end{cases}
            \label{eqn:ef-rom-1}\nonumber \\[0.3cm]
\end{eqnarray*}

\begin{eqnarray*}      
            &	\text{\bf (II)} &\text{\emph{Filter:}} \quad
\begin{cases} 
        	 -2 \delta^2 \, \Delta {\bu}_{n+1} +  {\bu}_{n+1} + \nabla \lambda_{n+1} = \bw_{n+1}& \text{in } \Omega \times \{t_{n+1}\}, \vspace{2mm}\\
             \nabla \cdot \bu_{n+1} = 0 & \text{in } \Omega \times \{t_{n+1}\}, \vspace{1mm}\\
{\bu}_{n+1} = \bu_D &\text{on } \Gamma_D \times \{t_{n+1}\}, \vspace{2mm}\\
\displaystyle \frac{\partial {\bu}_{n+1}}{\partial \boldsymbol n} = 0  & \text{on } \Gamma_N \times \{t_{n+1}\}.
\end{cases}
\end{eqnarray*}
The procedure is based on two steps: \textbf{(I)} solves the standard NSE for a given input $\bu_{n}$ at time $t_{n+1}$, providing the \textit{evolved velocity} $\bw_{n+1}$.
Step \textbf{(II)} applies a \textit{differential filter} (DF) of filter radius $\delta$, i.e., the spatial lengthscale of the filter action. 
The DF exploits an elliptic operator and smooths the small scales, i.e., high frequencies, of the intermediate flow velocity $\bw_{n+1}$, acting as a spatial filter. Finally, the \textit{filtered velocity} $\bu_{n+1}$ will be used as the next input in the time integration scheme. In the algorithm, $\lambda_{n+1}$ represents an auxiliary variable to guarantee the divergence-free constraint for $\bu_{n+1}$.
Following the notation of \cite{strazzullo2022consistency}, we denote with the acronym \textbf{noEF} the solution of the discretized NSE with no regularization.

\begin{remark}[Stokes filter]
\label{remStokesDF}
To preserve the incompressibility constraint during Step \textbf{(II)}, we rely on a Stokes DF \cite{ervin2012numerical, layton2008numerical}. Other filter actions might be used. For example, one could not enforce the divergence-free property $\nabla \cdot \bu_{n+1}$ and apply a div-grad stabilization \cite{heavner2017locally, john2017divergence}:
namely, in the DF equation, a term of the form
\begin{equation}
\tilde{\gamma} \nabla(\nabla \cdot {\bu}_{n+1}),
    \end{equation}
is added with $\tilde{\gamma}$ to be calibrated according to the test case. For an overview 
on this topic the reader may refer to \cite{decaria2018determination,john2017divergence,layton2009accuracy}. Despite the div-grad stabilization DF easier to code, in this contribution, we rely on the Stokes filter due to its consistency with respect to the NSE model and its robustness (no parameter calibration is needed except for $\delta$, see Remark \ref{rem:delta}). 
\end{remark}

\begin{remark}[Relaxation Step]
In the literature, another valuable algorithm has been employed to deal with convection-dominated settings: the Evolve-Filter-Relax (EFR) strategy. Indeed, in the EF approach, non-suited choices of the filter radius $\delta$ might yield overdiffusive simulations \cite{strazzullo2024variational,Ivagnes2025}. The EFR strategy follows the EF algorithm, but the new input $\bu_{n+1}$ is a convex combination of the evolved velocity field, $\bw_{n+1}$, and the filtered velocity field, say $\overline{\bw}_{n+1}$, i.e.,
$$
\bu_{n+1} = (1 - \chi)\bw_{n+1} + \chi \overline{\bw}_{n+1},
$$
for a chosen relaxation parameter $0 \leq \chi \leq 1$. However, in \cite{strazzullo2024variational} the authors question its applications due to spurious numerical oscillations arising in the relaxation step for a 2D flow past cylinder test case, i.e., one of the test cases we investigated in this contribution. Moreover, as already discussed in \cite{Ivagnes2025}, the parameter $\chi$ needs to be carefully chosen to avoid nonphysical simulations. Moreover, the optimization of both $\delta$ and $\chi$ in \cite{Ivagnes2025} gave comparable results to the optimization of only $\delta$ in the EF approach (see Remark \ref{rem:delta}). For these reasons, in this contribution, we propose to rely only on the EF algorithm to find the optimal $\delta$ using an RL-based strategy.  
\end{remark}
\begin{remark}[Classical choices for $\delta$]
    \label{rem:delta}
The \textit{optimal} value of $\delta$ is still an open question in CFD applications. In the literature, several scaling arguments have been proposed. For example, in \cite{strazzullo2022consistency,bertagna2016deconvolution}, the authors choose $\delta$ as the minimum element diameter $h_{min}$, or, for unstructured meshes, the Kolmogorov scale $\eta =L\cdot Re^{-\frac{3}{4}}$. 
In \cite{mou2023energy}, the choice is based on energy arguments.
\end{remark}
EF might suffer from a nontuned $\delta$ parameter: a too large value of $\delta$ yields overdiffusive results \cite{strazzullo2024variational}, while a too small value might not be sufficient to smooth out the spurious oscillations. 
In \cite{Ivagnes2025}, the authors propose the optimization of the filter radius in time, exploiting a few DNS data points as reference, to properly address the transition between laminar and turbulent regimes in numerical simulations. However, their data-driven strategy relies heavily on data, and no extrapolation results are proposed. 
In what follows, we propose an enhanced reference-guided data-driven strategy based on the use of RL capable \textit{of predicting the filter action in time} robustly.   
\section{Reinforcement Learning for Evolve-Filter}
\label{sec:rl-ef-all}
In this section, we provide a concise overview of Reinforcement Learning (RL) as a decision-making framework, introducing the Deep Q-Network (DQN) algorithm as an RL approach (Section \ref{sec:rl}), and describe how these concepts are applied to develop the RL-based Evolve--Filter (RL-EF) algorithm for optimized filtering in EF simulations (Section \ref{sec:rl-ef}).
From this section onward, $\bu_{n}$
will no longer denote the semi-discrete FEM approximation of the velocity field, but rather its vector whose entries are the coefficients of the FEM expansion in $\mathbb R^{N_{h}^{\bu}}$.

\subsection{An introduction to RL}
\label{sec:rl}

Reinforcement Learning (RL) is a paradigm based on the principle of \textit{learning by interaction}, where an agent learns to make decisions by interacting with a dynamic environment. The learning process is commonly formulated within the framework of a Markov Decision Process (MDP). More specifically, an MDP is defined by a tuple $(\mathcal{S}, \mathcal{A}, P, r, \gamma)$, where $\mathcal{S}$ denotes the state space, $\mathcal{A}$ the action space, $P$ the state transition probability, $r$ the reward function, and $\gamma \in (0,1]$ the discount factor. The interaction between agent and environment is organized in \textit{episodes}, each consisting of a finite sequence of training decision \textit{steps}. 
At a given decision step $t_n$, the environment is characterized by a state $s_n = s(t_n) \in \mathcal{S}$. Based on the observation of the current state, the agent selects an action $a_n = a(t_n) \in \mathcal{A}$ according to a policy $\pi(a_n \mid s_n)$, which defines the probability of choosing the action $a_n$ given the state $s_n$. As a consequence of the selected action, the environment transitions to a new state $s_{n+1}$ following the transition law $P(s_{n+1} \mid s_n, a_n)$, and the agent receives a scalar reward $r_{n+1}$.
The objective of the agent is to learn an optimal policy that maximizes the expected cumulative (possibly discounted) reward over time. In particular, for a given episode, the cumulative reward (or return) is defined as
\[
R_{\text{train}}= \sum_{n=0}^{N_{\text{train}}-1} \gamma^{n} r_{n+1},
\]
where $N_{\text{train}}$ denotes the episode length, i.e., the number of time steps in each episode, and $\gamma \in (0,1]$ is the discount factor. Through repeated interactions with the environment, the agent progressively improves its policy, enabling generalization to states that may not have been explicitly encountered during training, and, hence, time extrapolation capability.

\subsubsection{The Deep-$Q$-Network algorithm}
In this contribution, we adopt Deep $Q$-Network (DQN)~\cite{mnih2015human}, a value-based RL algorithm which models the $Q$-function, i.e., the action--value function:
\[
Q^\pi(s,a) = \mathbb{E}_\pi \left[ \sum_{k=0}^{\infty} \gamma^k r_{n+k} \,\middle|\, s_n=s, a_n=a \right].
\]
In the above expression, $\mathbb{E}_\pi[\cdot]$ denotes the expectation with respect to the stochasticity induced by the policy $\pi$ and the environment dynamics, i.e., over all possible future state--action trajectories generated by following policy $\pi$ starting from $(s,a)$. The infinite sum represents the discounted cumulative reward collected along the future trajectory, where the discount factor $\gamma$ ensures convergence and assigns decreasing importance to rewards received at later time steps.
In DQN, the $Q$-function is approximated by means of a deep neural network. The $Q$-Network represents a nonlinear mapping
\[
Q_\theta : \mathcal{S} \times \mathcal{A} \rightarrow \mathbb{R},
\]
which estimates the expected cumulative reward associated with each state--action pair.
The policy is not learned explicitly, but it is instead induced from the $Q$-function through a greedy selection rule,
\[
a_n = \pi(s_n) = \arg\max_{a \in \mathcal{A}} Q_\theta(s_n,a),
\]
possibly combined with an exploration strategy during training.
In DQN, two neural networks are employed: an online $Q$-Network with parameters $\theta$, which is updated during training, and a target $Q$-Network with parameters $\theta^-$, which is periodically synchronized with the online network to make the learning process more stable.
The $Q$-Network parameters $\theta$ are optimized by minimizing the temporal-difference loss
\[
\mathcal{L}(\theta) =
\mathbb{E} \Big[
\big( r_{n+1} + \gamma \max_{a'} Q_{\theta^-}(s_{n+1},a') - Q_\theta(s_n,a_n) \big)^2
\Big],
\]
where $\theta^-$ denotes the parameters of the target network.
For a deeper insight into RL, the reader may refer to \cite{RL}.

\begin{figure}[htpb!]
    \centering
    \includegraphics[width=\linewidth, trim={4cm 4cm 4cm 4cm}, clip]{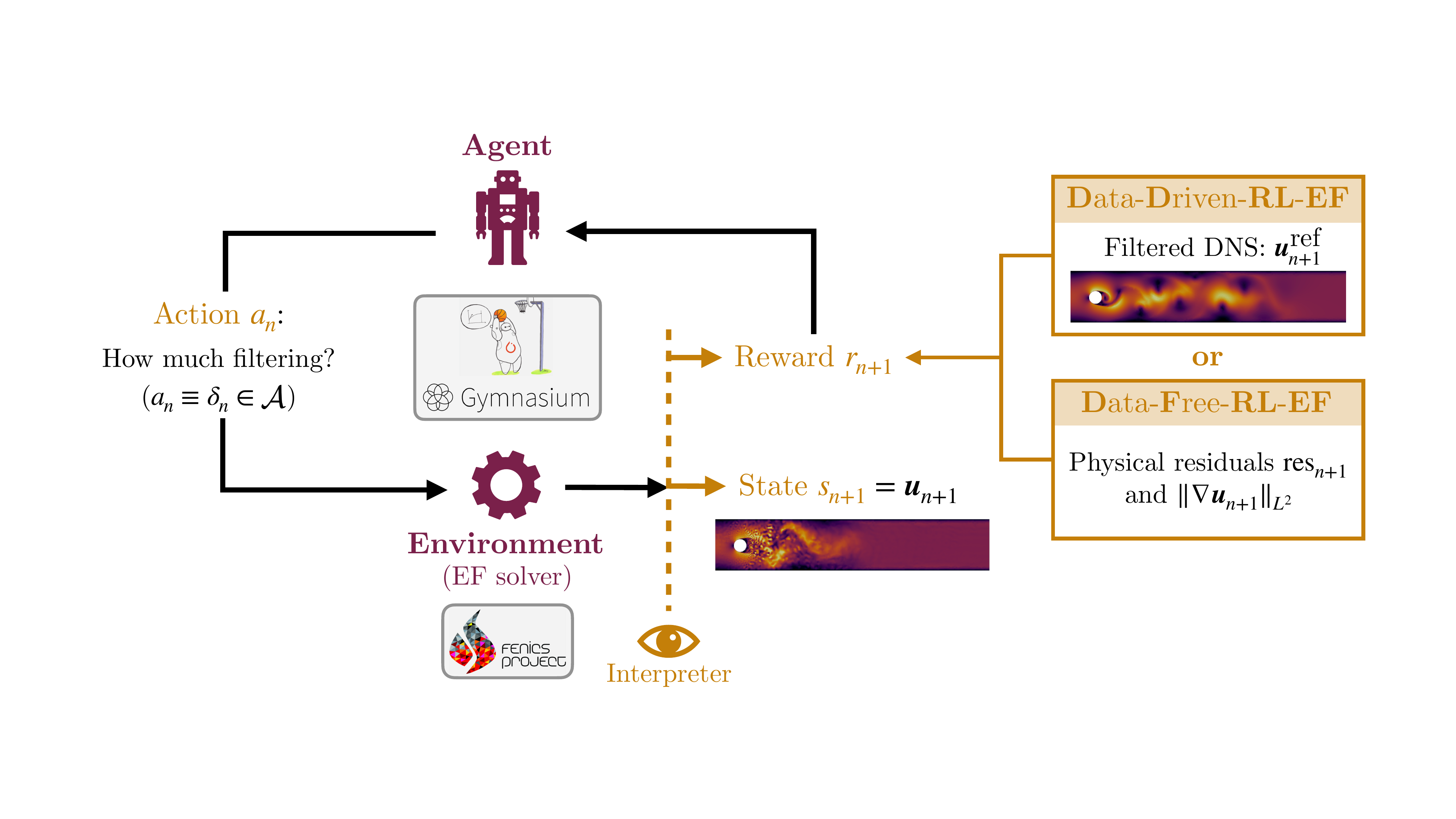}
    \caption{Schematic time stepping mechanism in an RL-EF training episode, and related software.}
    \label{fig:RLpipe}
\end{figure}

\subsection{The RL-EF algorithm}
\label{sec:rl-ef}
\begin{figure}[htpb!]
    \centering
    \includegraphics[width=0.8\linewidth]{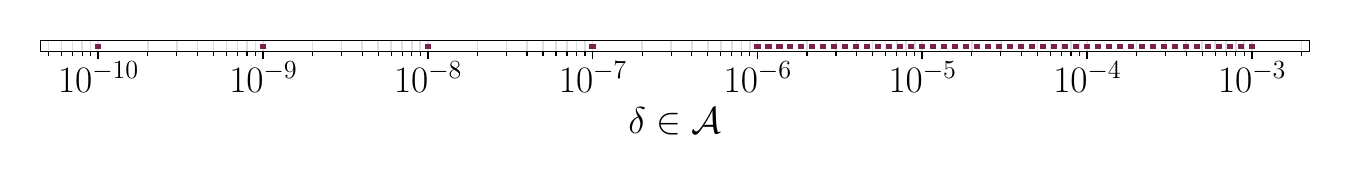}
    \caption{Action space in RL-EF.}
    \label{fig:action-space}
\end{figure}

We apply the RL techniques described in Section \ref{sec:rl} in the context of EF simulations, introducing the Reinforcement Learning-based Evolve-Filter algorithm (\textbf{RL-EF}).
More in detail, RL-EF is an optimized EF simulation where the filter radius $\delta_n$ is time-dependent and updated at each time step by a pre-trained RL agent.
In the proposed setting, the action coincides with the filter radius to be applied at time $t_n$, i.e., $a_n=\delta_n$, while the observation (state) $s_{n+1}$ corresponds to the vector-valued filtered velocity field $\boldsymbol{u}_{n+1} (\delta_n)$ obtained using the filter radius selected by the RL agent.
Hence, the action space $\mathcal{A}$ takes $50$ discrete values, composed of four logarithmically-spaced values in $[10^{-10}, 10^{-7}]$ and fortysix logarithmically-spaced values in $[10^{-6}, 10^{-3}]$, schematically represented in Figure \ref{fig:action-space}. The logic behind this choice is the effect of the filter on the velocity field: small filters values with orders of magnitudes from $\mathcal{O}(10^{-10})$ to $\mathcal{O}(10^{-7}$) produce similar EF solutions, e.g., $\bu_{n}(10^{-10})\approx \bu_{n}(10^{-7})$, while larger values have a significant influence on the solution, e.g. $\bu_{n}(10^{-6}) \napprox \bu_{n}(10^{-3})$. For this reason, we consider more values in the range $[10^{-6}, 10^{-3}]$.

The observation space depends instead on the number of degrees of freedom of the discretization adopted in the EF simulation, namely $\mathcal{S} \equiv \mathbb{R}^{N_h^{\bu}}$.
In our setting, an RL episode corresponds to a sequence of physical time steps, i.e., the time steps of an EF simulation over a prescribed training time window, as represented in Figure \ref{fig:RLpipe}.
As a result of the training procedure, the RL agent identifies an optimized trajectory
\[
\{ (\boldsymbol{u}_{1}, \delta_{1}), \dots, (\boldsymbol{u}_{N}, \delta_{N}) \}
\]
over the time interval $(0,T]$, by learning from a reduced training window $(0,T_{\text{train}}]$, with $T_{\text{train}} \ll T$.
\begin{remark}[Proximal Policy Optimization]
Other RL strategies might be employed, such as the Proximal Policy Optimization (PPO)
algorithm \cite{schulman2017ppo}, which allows us to consider a continuous action space. We also tried this approach in our test cases; however, we obtained poor results in terms of solution accuracy and training stability. We do not exclude that a
{more advanced empirical investigation and fine-tuning} of the net architecture might yield more accurate results. However, this goes beyond the main scope of the paper.
\end{remark}

\subsubsection{RL-EF rewards}
A natural question arises in the context of RL: \textit{what is the most appropriate reward to let the system learn a suitable filter action?}
In our numerical investigation, we consider two types of reward: either \textbf{data-driven}, based on the accuracy of the RL-EF solution with respect to some pre-collected reference data, or \textbf{data-free}, based on the physical residual and on the RL-EF velocity gradients. Consequently, the general RL-EF methodology specializes in \textbf{DD-RL-EF} (Data-Driven Reinforcement-Learning-based Evolve-Filter) and \textbf{DF-RL-EF} (Data-Free Reinforcement-Learning-based Evolve-Filter).
In the data-driven setting, we use the data from a reference DNS simulation, which is performed on a refined grid, and then projected into the \textit{marginally-resolved} grid used in RL-EF.
From now on, we will refer to the reference velocity solution as \textbf{Filtered DNS}, and we will indicate it as $\boldsymbol{u}_n^{\text{ref}}$ at time step $t_n$.
\medskip

In terms of reward, as in~\cite{kurz2022deep, kurz2023deep, kurz2025harnessing}, we consider the following general form for the reward at each time step $t_{n+1}$:
\begin{equation}
    r^{\text{DD}}_{n+1} = 2 \exp{\left( -\dfrac{1}{\alpha} \mathrm{e}_{n+1}\right)} - 1,
    \label{eq:general-form-reward}
\end{equation}
where $\mathrm{e}_{n+1}$ is a metric associated with the error of the RL-EF solution $\bu_{n+1}(\delta_n)$ with respect to the DNS reference $\bu^{\text{ref}}_{n+1}$. The expression \eqref{eq:general-form-reward} is widely used in RL to map the error $\mathrm{e}_{n+1}$ into the interval $[-1, 1]$. This transformation improves stability and facilitates convergence during the training stage.
\begin{remark}
In the literature, when applying RL to fluid dynamics applications, the error metric taken into account is typically the relative energy spectrum difference at the first (and more energetic) wave numbers, namely:
\begin{equation}
    \mathrm{e}^{\text{spectrum}}_{n+1} = \dfrac{1}{\kappa_{\text{max}}} \sum_{\kappa=1}^{\kappa_{\text{max}}} \left( \dfrac{\mathcal{E}_{n+1}(\kappa)-\mathcal{E}_{n+1}^{\text{ref}}(\kappa)}{\mathcal{E}_{n+1}^{\text{ref}}(\kappa)} \right)^2,
    \label{eq:err-spectrum}
\end{equation}
where $\mathcal{E}(\kappa)$ and $\mathcal{E}^{\text{ref}}(\kappa)$ are the RL-EF and the filtered DNS energies, respectively, at the wavenumber $\kappa$ at time step $t_{n+1}$. The value $\kappa_{\text{max}}$ is the maximum wavenumber considered in the error, which is usually selected via empirical tests and depends on the test case.
Preliminary experiments demonstrate that the spectrum error metric \eqref{eq:err-spectrum} does not provide good results for RL-EF, in our test cases.
Hence, we investigate the following data-driven error, which coincides with the relative error between the RL-EF velocity and the reference solution:
\begin{equation}
    \mathrm{e}^{\text{DD}}_{n+1}=\left(\dfrac{\|\bu_{n+1} - \bu_{n+1}^{\text{ref}}\|_2}{\|\bu_{n+1}^{\text{ref}}\|_2}\right)^2,
\end{equation}
where the notation $\|(\cdot)\|_2$ stands for the $L^2$ norm in the domain $\Omega$.
\end{remark}
In our tests, the \textbf{DD-RL-EF} reward is $r_{n+1}^{\text{DD}}$ in Equation \eqref{eq:general-form-reward}, with $\alpha=1$ and $\mathrm{e}_{n+1}=\mathrm{e}^{\text{DD}}_{n+1}$.
Two questions automatically arise for DD-RL-EF:

\begin{itemize}
    \item[(i)] \textit{What is the number of time steps that should be considered in one RL episode?}
    \item[(ii)] \textit{Should each episode start from the beginning of the simulation, e.g. $t=0$?}
\end{itemize}
To answer these questions, we considered two variants of DD-RL-EF:
\begin{itemize}
    \item a standard DD-RL-EF, which we refer to simply as \textbf{DD-RL-EF}, where each episode starts from $t=0$ and has duration $T_{\text{train}}$ (corresponding to $N_{\text{train}}$ time steps). In our experiments, we set $T_{\text{train}}=T/4$. In this case, the cumulative reward for each episode will be:

    \begin{equation}
        R^{\text{DD}}=\sum_{n=1}^{N_{\text{train}}} r_{n}^{\text{DD}}.
        \label{eq:cum-rew-dd-rl}
    \end{equation}
    
    \item \textbf{DD-RL-EF-rand}, where each episode starts from a randomly selected time instant within the training window $n_{\text{rand}} \in \{ 1, 2, \dots,  , N_{\text{train}} - 1, N_{\text{train}}\}$. In this setting, the episode length is reduced compared to the standard DD-RL-EF case, and it is set to $T_{\text{train}}^{\text{rand}}=T/10$, corresponding to $N^{\text{rand}}_{\text{train}}$ time steps. This strategy promotes exploration of a wider range of states than standard episodes initialized at $t=0$. The cumulative reward would be:

    \begin{equation}
        R^{\text{DD}}_{\text{rand}}=\sum_{n=n_{\text{rand}}}^{N^{\text{rand}}_{\text{train}}} r_{n}^{\text{DD}}.
        \label{eq:cum-rew-dd-rl-rand}
    \end{equation}
\end{itemize}
\medskip
On the other hand, in the \textit{data-free} setting, we consider the following reward expression:

\begin{equation}
    r_{n+1}^{\text{DF}} = \dfrac{1}{2} \Bigl( 2\exp{\left( -\alpha_{\text{res}} \, \mathrm{res}_{n} \right) -1}\Bigr) + \dfrac{1}{2} \left(  2 \exp{\left( - \left| \dfrac{1}{\alpha_{\text{grad}} \|\nabla \bu_{n+1}\|_2 - \|\nabla \bu_{n}\|_2)} \right| \right) -1}\right),
\end{equation}
which includes two contributions: 
$$
\text{(i) } 2\exp{\left( -\alpha_{\text{res}} \, \mathrm{res}_{n+1} \right) -1}, \text{ and  (ii) } 2 \exp{\left( - \left| \dfrac{1}{\alpha_{\text{grad}} (\|\nabla \bu_{n+1}\|_2 - \|\nabla \bu_{n}\|_2)} \right| \right) -1},
$$
for $\alpha_{\text{res}}$ and $\alpha_{\text{grad}}$ positive real numbers. Both parts share the same form as given in Equation \eqref{eq:general-form-reward}:

\begin{itemize}
    \item[(i)] The contribution is directly related to the residual of the FE system $\mathrm{res}_{n+1}$. This leads to decreasing rewards as the physical residual increases, i.e., when the RL-EF prediction is not a solution of the discrete NSE in Equation \eqref{eq:NSE}. This component penalizes excessively small values of $\delta_n$, which do not stabilize the system. As observed in the noEF simulation, such behaviour leads to under-resolved equations and a subsequent blow-up of the residual.
    \item[(ii)] The contribution is instead related to the temporal change in the velocity space gradients. Small changes of the gradients manifest in overdiffusive behaviour, as in the EF simulation with a large $\delta$ value, such as $\delta=\eta$. Unlike the residual term, this component penalizes too large $\delta_n$ values. The agent indeed gets a larger reward when the velocity gradients change more over time, better capturing the flow features expected by its variation in time. 
\end{itemize}
The constants $\alpha_{\text{res}}$ and $\alpha_{\text{grad}}$ are tuned after preliminary empirical investigations.
In our tests, the \textbf{DF-RL-EF} reward is $r^{\text{DF}}_{n+1}$, with $\alpha_{\text{res}}=\num{1e5}$ and $\alpha_{\text{grad}}=\num{1e4}$. Hence, we indicate the cumulative reward for each episode as:
\begin{equation}
    R^{\text{DF}}=\sum_{n=1}^{N_{\text{train}}}r_{n}^{\text{DF}}.
\end{equation}

\begin{remark}[Why not DF-RL-EF-rand?]
    To preserve a fully data-free approach, we do not adopt the strategy of initializing episodes at a random time step $n=n_{\text{rand}}$, as this would require access to the filtered DNS state $\bu_{n_{\text{rand}}}^{\text{ref}}$ to start the time integration without compromising the accuracy of the prediction.
\end{remark}
Finally, the reward formulation can be adapted to incorporate weak structure-preserving constraints, like energy or enstrophy conservation, which do not require access to the reference data, but only to the RL-EF solution.
In this context, we consider a Structure-Preserving variant of DF-RL-EF (\textbf{SP-DF-RL-EF}), where the reward includes some penalization terms for the energy ($\mathcal E_n$ at time $t_n$) and enstrophy ($\mathcal Z_n$ at time $t_n$) growth:

\begin{equation*}
    R^{\text{DF}}_{\text{SP}}=\sum_{n=1}^{N_{\text{train}}} \left[ r_{n}^{\text{DF}} + \mathbf{1}_{\{\mathcal{E}_{n}>\mathcal{E}_{n-1}\}}\left( \exp{\left( -\dfrac{\mathcal{E}_{n} - \mathcal{E}_{n-1}}{\alpha_{\mathcal{E}}\, \mathcal{E}_{n-1}} \right)} -1 \right) + \mathbf{1}_{\{\mathcal{Z}_{n}>\mathcal{Z}_{n-1}\}}\left( \exp{\left( -\dfrac{\mathcal{Z}_{n} - \mathcal{Z}_{n-1}}{\alpha_{\mathcal{Z}}\, \mathcal{Z}_{n-1}} \right)} -1 \right) \right].
    \label{eq:sp-reward}
\end{equation*}
In the above expression, $\mathbf{1}_{\{\mathcal{E}_{n}>\mathcal{E}_{n-1}\}}$ and $\mathbf{1}_{\{\mathcal{Z}_{n}>\mathcal{Z}_{n-1}\}}$ indicate the indicator functions that are activated only when there is an energy and an enstrophy growth, respectively. In such cases, an additional penalization term related to the energy and/or enstrophy growth is added into the reward. We tuned the parameters associated with the penalization terms to $\alpha_{\mathcal{E}}=\alpha_{\mathcal{Z}}=0.1$.
Table \ref{tab:acronyms-all} summarizes the acronyms used for the RL-EF methodologies introduced in this work.

\begin{table*}[h]
\centering
\caption{Acronyms and main features of the proposed RL-EF approaches.}
\label{tab:acronyms-all}
{
\begin{tabular}{ccccc}
\toprule
{\textbf{Acronym}} & \textbf{Episode start}& \textbf{Episode length} &\textbf{Cum. reward}&\textbf{Data needed}
\\ \midrule
DD-RL-EF & $t=0$ & $T_{\text{train}}=T/4$& $R^{\text{DD}}$ & $\bu^{\text{ref}}_n, t_n \in [0, T_{\text{train}}]$\\  
DD-RL-EF-rand & random $t \in [0, T_{\text{train}}]$ & $T_{\text{train}}^{\text{rand}}=T/10$& $R^{\text{DD}}_{\text{rand}}$ & $\bu^{\text{ref}}_n, t_n \in [0, T_{\text{train}} + T_{\text{train}}^{\text{rand}}]$\\  
DF-RL-EF & $t=0$ & $T_{\text{train}}=T/4$& $R^{\text{DF}}$ & No reference needed\\ 
SP-DF-RL-EF& $t=0$ & $T_{\text{train}}=T/4$& $R^{\text{DF}}_{\text{SP}}$ & No reference needed\\
\bottomrule
\end{tabular}}
\end{table*}

\paragraph{\textbf{Software}: }
As indicated in Figure \ref{fig:RLpipe}, the libraries used to develop RL-EF are \textbf{Gymnasium}~\cite{towers2024gymnasium, kwiatkowski2024gymnasium}, a standardized Python-based library for defining reinforcement learning environments, and \textbf{Stable Baselines 3}~\cite{raffin2021stable}, a library of reliable, well-tested reinforcement learning algorithms built on top of Gymnasium. The RL environment is integrated with the open-source Python-based FEM framework FEniCSX~\cite{baratta2023dolfinx, alnaes2014unified, scroggs2022construction, scroggs2022basix}, which is used to define the simulation setup and the EF time stepping.

\section{Numerical results}
\label{sec:results}
This section is dedicated to highlighting the potential of the RL-EF methodologies in the following two test cases:
\begin{itemize}
    \item[(i)] two-dimensional flow past a cylinder at Reynolds number $Re=\num{1000}$ (Section \ref{sec:test-case-1});
    \item[(ii)] two-dimensional decaying homogeneous turbulence test case, at $Re=\num{40000}$ (Section \ref{sec:test-case-2}).
\end{itemize}
We recall here the notation used for the different types of simulation:
\begin{itemize}
    \item \textbf{DNS}: expensive simulation on the most refined grid.
    \item \textbf{Filtered DNS}: reference data used for RL reference-guided training and for comparison, obtained by interpolating the DNS onto the function space defined on the coarse mesh.
    \item \textbf{noEF}: DNS on the coarse mesh, which, depending on the grid resolution, might blow up and needs to be stabilized.
    \item \textbf{EF (Kolmogorov)}: the standard EF approach applied to the coarse mesh, where $\delta_n = \delta=\eta \quad \forall n=1,\dots, N$.
    \item \textbf{RL-EF}: novel EF simulations on the coarse mesh in which a small time-window pre-trained RL agent predicts the time-dependent filter radius $\delta_n$ for the unseen time instances. The acronyms and the corresponding features of the data-driven and data-free RL-EF strategies have been introduced in Table \ref{tab:acronyms-all}. 
\end{itemize}

\subsection{Test case 1: flow past a cylinder}
\label{sec:test-case-1}
\begin{figure}[htpb!]
    \centering
    \subfloat[Coarse mesh]{\includegraphics[width=0.5\linewidth, trim={0 10cm 0 10cm}, clip]{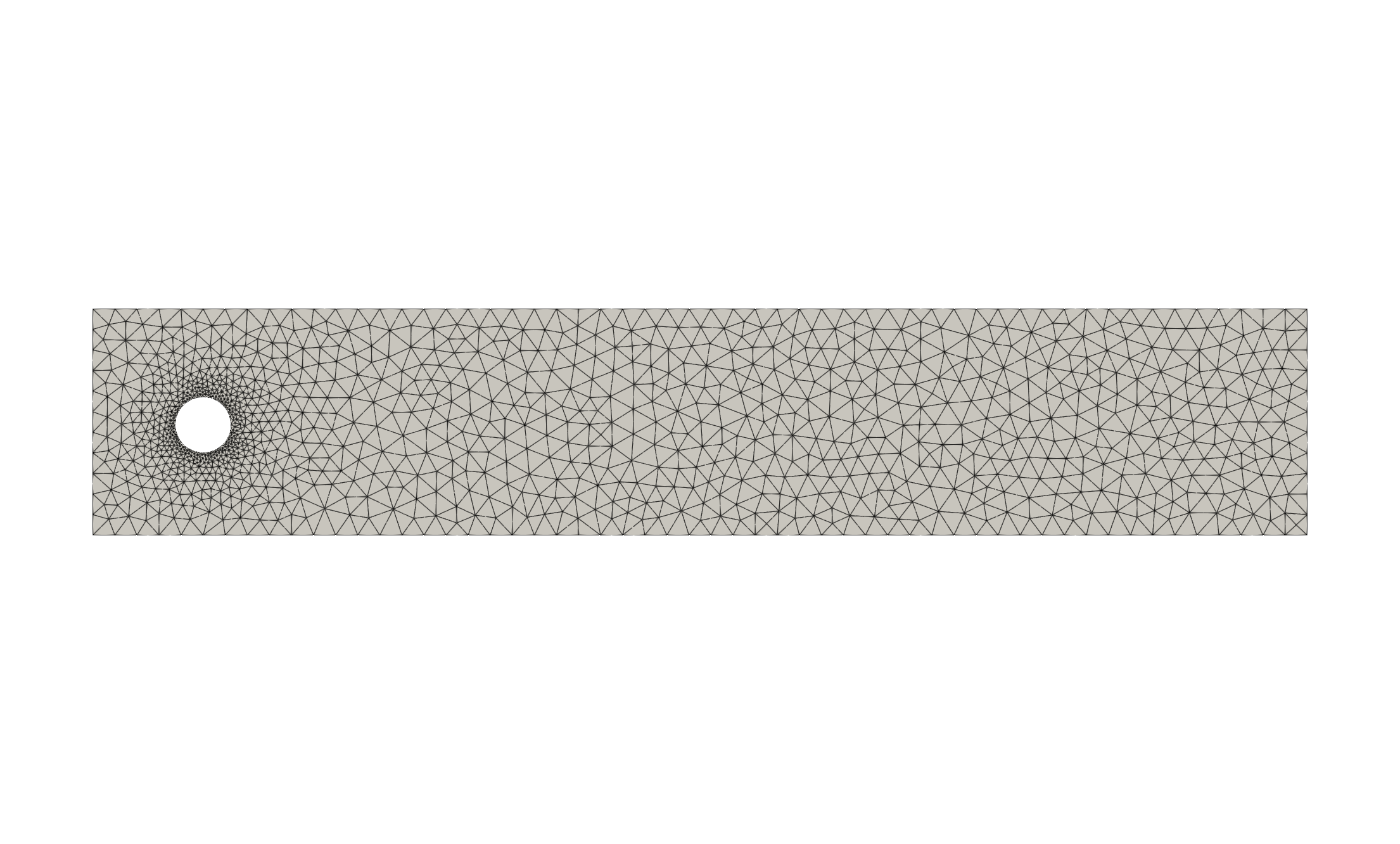}}
    \subfloat[Fine mesh]{\includegraphics[width=0.5\linewidth, trim={0 10cm 0 10cm}, clip]{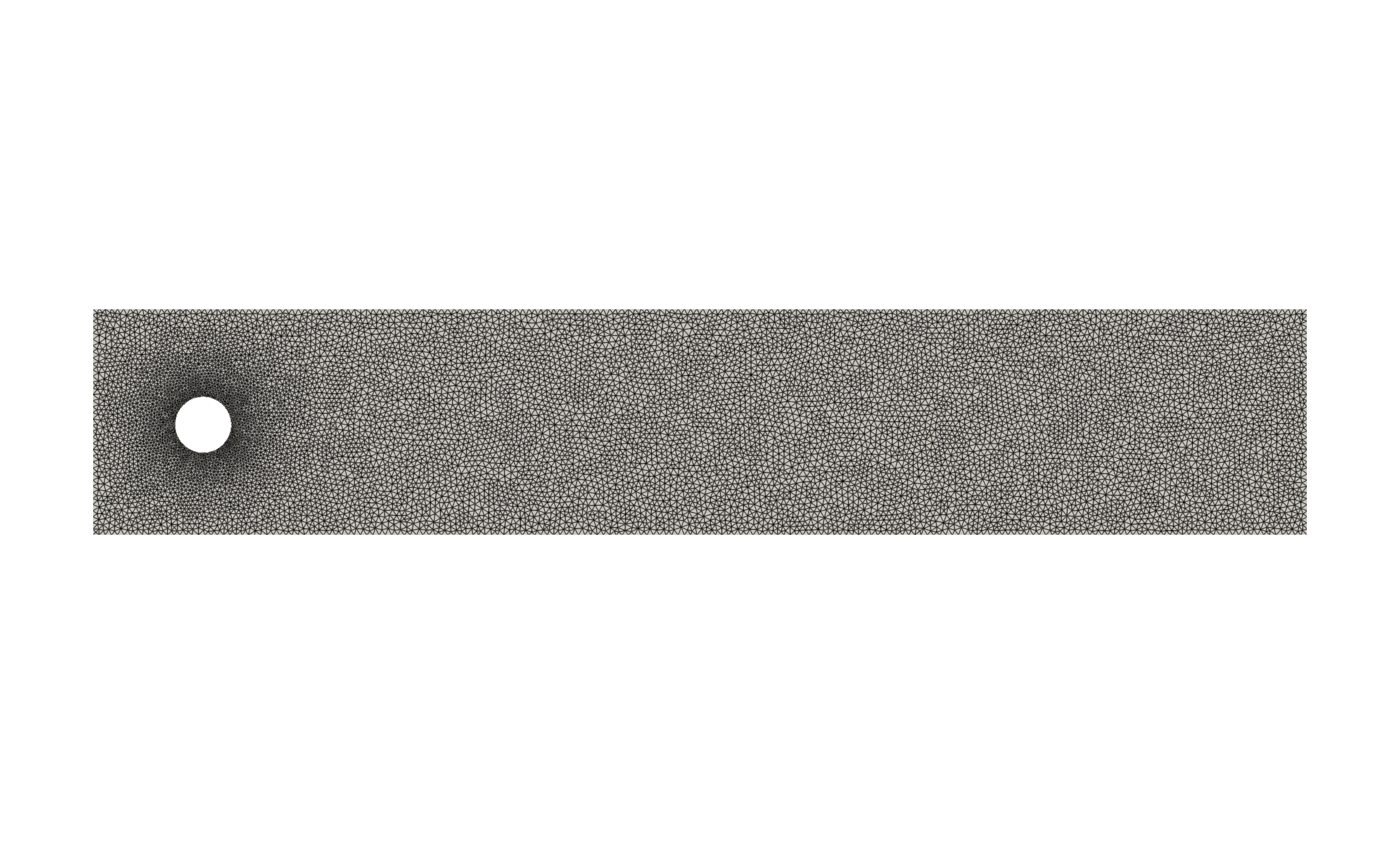}}
    \caption{\emph{Test case 1}. The two grids used for the simulations. The coarse grid is used for noEF, EF, and RL-EF simulations, while the fine mesh is used for the DNS reference simulation.}
    \label{fig:meshes-cyl}
\end{figure}

We consider the test case of the incompressible flow past a circular cylinder at $Re=1000$, constant in time.
Table \ref{tab:dofs-cyl} reports the number of degrees of freedom used for the velocity ($N_h^{\boldsymbol{u}}$)
and for the pressure ($N_h^p$) field, and the minimum and maximum mesh sizes, for both the fine grid used for the DNS simulation, and the marginally-resolved grid used in noEF, EF, RL-EF approaches. The grids are also displayed in Figure \ref{fig:meshes-cyl}, and are characterized by triangular elements.
The simulations are performed in the time window $[0, T]$, with $T=4$.
The time step used for simulations on the coarse grid is $\Delta t = \num{4e-4}$, for a total of $N=\num{10000}$ time instances.
\begin{table}[H]
    \caption{\emph{Test case 1}. Number of degrees of freedom, minimum and maximum sizes for the grids.}
    \centering
    \begin{tabular}{ccccc}
    \toprule
    & $N_h^{\boldsymbol{u}}$ & $N_h^p$ & $h_{min}$ & $h_{max}$ \\
    \midrule
        Coarse mesh & $9792$ & $1274$& $4.46e-3$&$4.02e-2$  \\
         Fine mesh & $43609$ & $11036$&$4.44e-3$&$1.13e-2$ \\
    \bottomrule     
    \end{tabular}
    \label{tab:dofs-cyl}
\end{table}

The domain is $\Omega = \{(x, y) \in [0, 2.2] \times [0, 0.41] \text{ such that } (x-0.2)^2 +(y-0.2)^2 \geq 0.05^2\}$. We consider no-slip boundary conditions on $\partial \Omega_{wall}$ (which includes the top, bottom, and cylinder boundaries), and an inlet velocity $\mathbf{u}_{in}$ (left boundary), given by:
\begin{equation}
    \bu_{in }= \left( \dfrac{6}{0.41^2}\,y\,(0.41-y), 0
    \right)\,.
    \label{u_inlet}
\end{equation}

\begin{figure}[htpb!]
    \centering
    \subfloat{\includegraphics[width=.8\linewidth, trim={15cm 0 0 0}, clip]{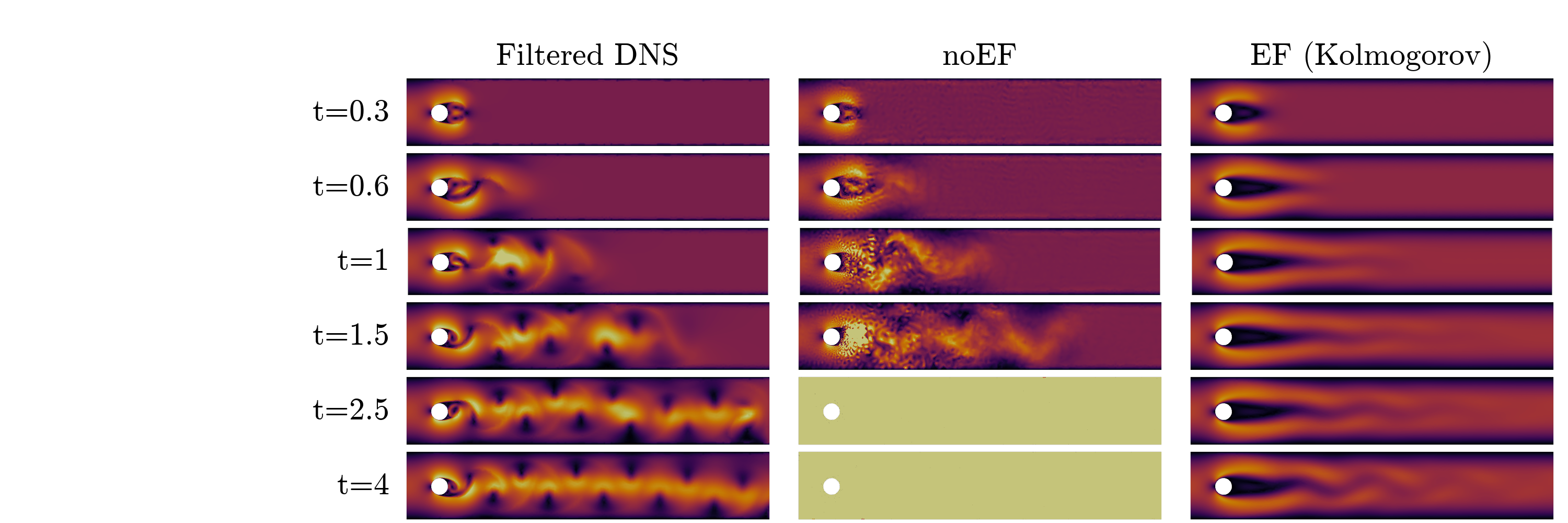}}
    \subfloat{\includegraphics[width=.2\linewidth, trim={90cm 31cm 25cm 0cm}, clip]{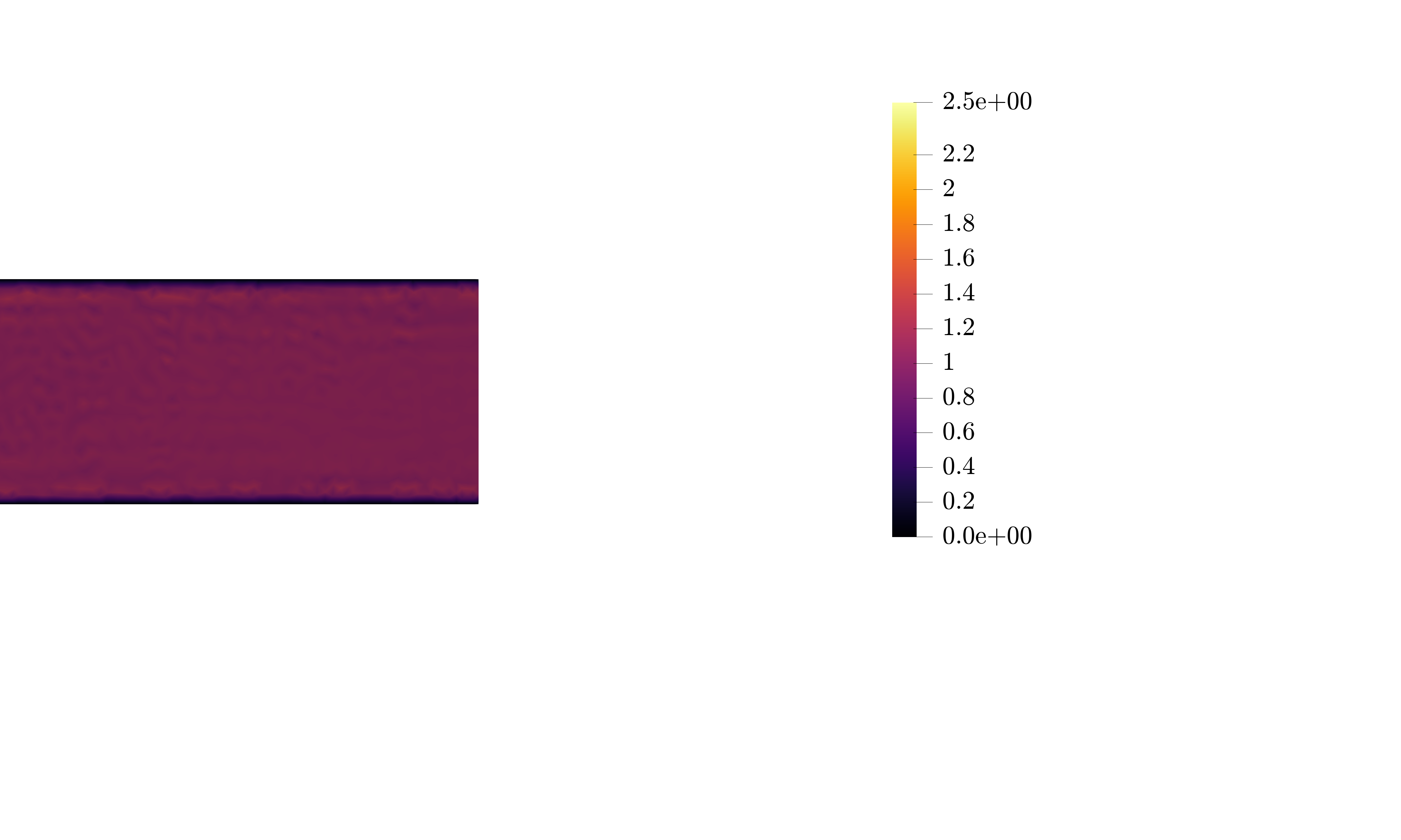}}
    \caption{\emph{Test case 1}. Velocity magnitude fields for the filtered DNS, the noEF and EF with $\delta=\eta$. The fields are represented at different time instances.}
    \label{fig:cylinder-pre-graphics}
\end{figure}

On the outlet of the domain (right boundary), we employ standard free-flow conditions.  The initial condition is $\bu_0=(0, 0)$. The viscosity is set to $\nu=\num{1e-4}$.
Figure \ref{fig:cylinder-pre-graphics} shows representative solutions in time of the velocity magnitude for the filtered DNS, the noEF and EF (Kolmogorov) cases, and it demonstrates that:
\begin{itemize}
    \item noEF exhibits noisy behaviour, and after $t=1.5$ the energy grows unbounded and the simulation blows up;
    \item EF (Kolmogorov) has an overdiffusive behaviour, as it manifests a vortex shedding in the cylinder wake only at the end of the simulation, after $t=2$.
\end{itemize}

\begin{figure}[htpb!]
    \centering
    \includegraphics[width=\linewidth]{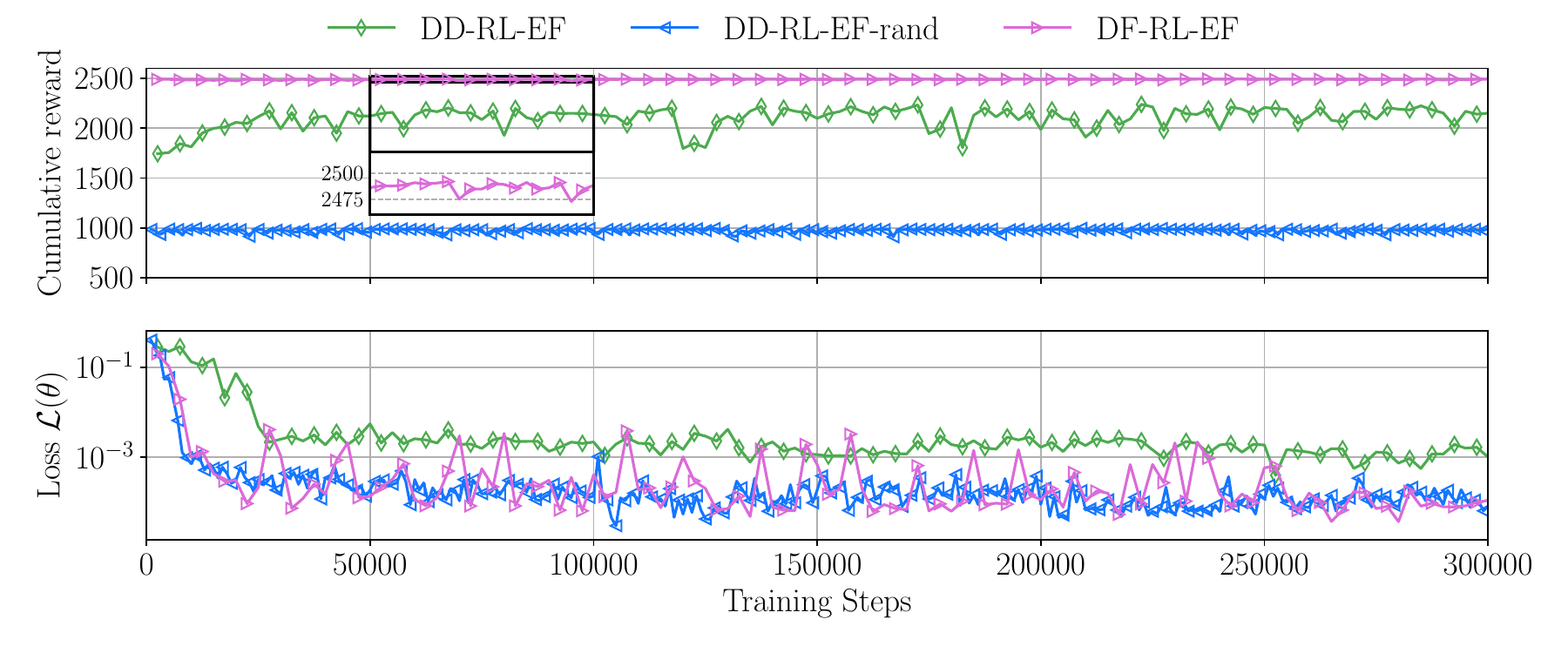}
    \caption{\emph{Test case 1}. Cumulative reward and loss value during training steps, for the RL-EF methods.}
    \label{fig:cylinder-rewards}
\end{figure}

The poor result obtained with noEF and EF (Kolmogorov) motivates the use of RL-EF.
Before showing the qualitative results of RL-EF simulations, we analyze the RL training statistics in Figure \ref{fig:cylinder-rewards}, and in particular, the cumulative reward over the episodes and the DQN loss function over the time steps.
First of all, it is important to highlight that the maximum cumulative reward that can be reached coincides with the number of time steps for each episode $N_{\text{train}}$. In DD-RL-EF and DF-RL-EF, $N_{\text{train}}=N/4 = 2500$, corresponding to $T_{\text{train}}=1$, while in DD-RL-EF-rand, as specified in section \ref{sec:rl-ef}, $N_{\text{train}}^{\text{rand}}=1000$, corresponding to $T_{\text{train}}^{\text{rand}}=0.4$.
It can be immediately noted that all the methods converge close to their maximum achievable cumulative reward, indicating that they successfully learn near-optimal policies.
DF-RL-EF reaches the performance plateau earlier and exhibits smaller fluctuations compared to DD-RL-EF, as can be seen from the zoomed box in the rewards of Figure \ref{fig:cylinder-rewards}. This depends on the reward expression used in DF-RL-EF, which includes the residual contribution and the time change of the velocity gradients. This form, together with the choice of the parameters $\alpha_{\text{res}}$ and $\alpha_{\text{grad}}$, is specifically designed to obtain a stable RL learning dynamics.
Also in DD-RL-EF-rand training, the cumulative reward exhibits small oscillations and quickly reaches its maximum value ($1000$). The random selection of the initial time for each episode yields a larger state exploration and accelerates the agent's learning.
This discussion is confirmed by the loss evolution across episodes and time steps: the DD-RL-EF-rand and DF-RL-EF losses decay faster than DD-RL-EF, and reach $\sim \num{1e-4}$ before episode $100$. On the other hand, the DD-RL-EF loss reaches a larger plateau of $\num{1e-3}$.

\begin{figure}[htpb!]
    \centering
    \includegraphics[width=\linewidth]{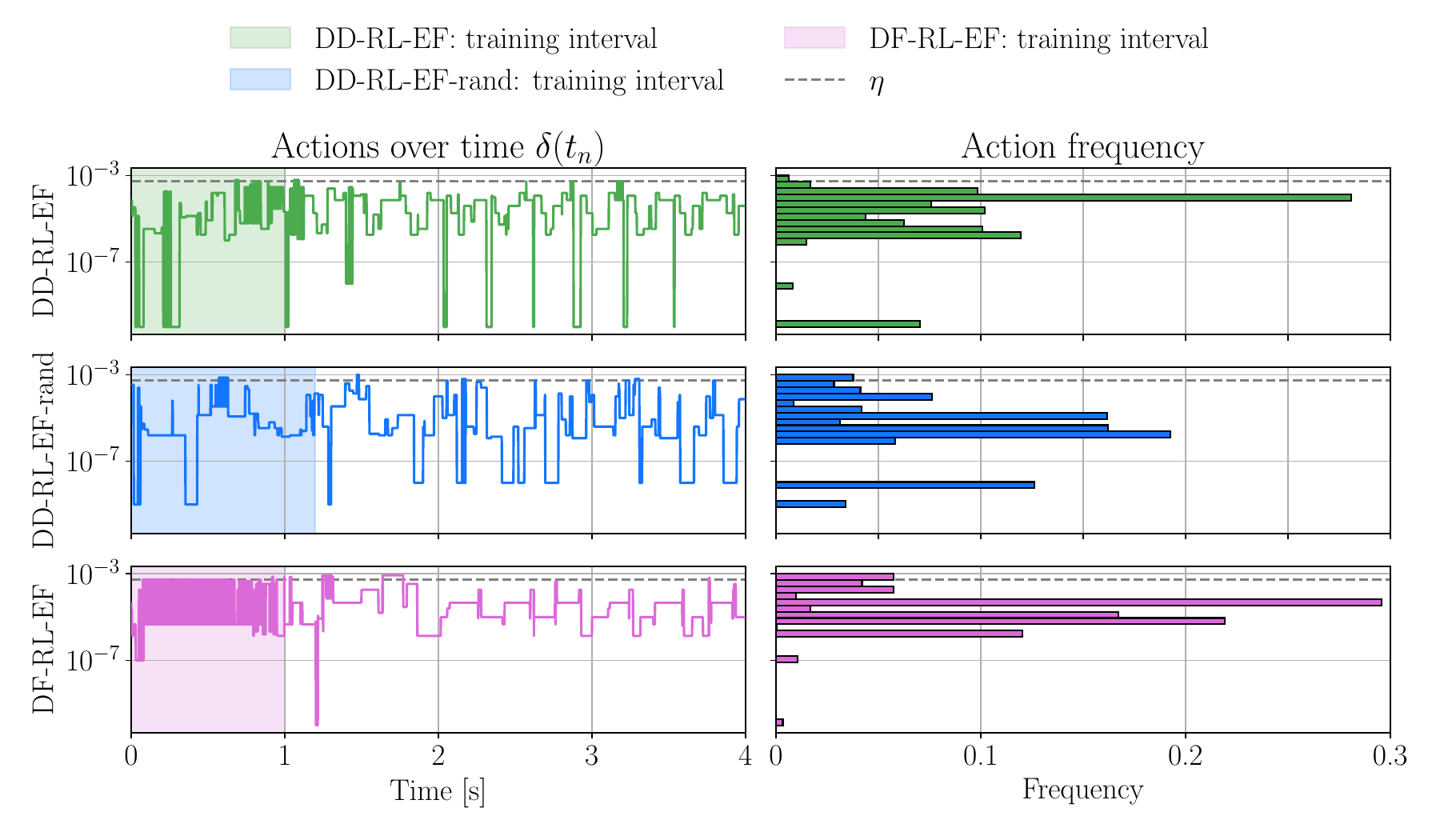}
    \caption{\emph{Test case 1}. Action values over testing time, and corresponding frequency, for the RL-EF methods.}
    \label{fig:cylinder-actions}
\end{figure}

Once completed the RL training, we can analyse the performance of RL-EF strategies on a \textit{online} simulation, with time extrapolation. Figure \ref{fig:cylinder-actions} displays: (i) the actions $\delta_n, n=1, \dots, N$ selected by the pre-trained RL agent over the simulation time (on the left); (ii) the action frequency, answering the question: \textit{what is the percentage of time an action is chosen by the agent?}

First, we can immediately see that different agents choose different actions, but in all algorithms, the action values have an oscillatory and non-regular time behaviour. Similar to the preliminary investigation conducted in~\cite{Ivagnes2025}, this suggests that the simulation does not have a unique and fixed optimal $\delta$ for the entire time window, and that, depending on the simulation state, for some time instances the optimal choice is \textit{not to filter}, i.e., $\delta_n$ of order of magnitude from $\mathcal{O}(10^{-10})$ to  $\mathcal{O}(10^{-7})$. The RL agent choice is indeed most of the time a smaller value than the Kolmogorov length scale $\eta$, indicated in Figure \ref{fig:cylinder-actions} with dashed grey lines.
Moreover, DF-RL-EF converges overall to larger values of $\delta_n$ compared to DD-RL-EF and DD-RL-EF-rand. This behaviour suggests that the data-free reward formulation implicitly assigns greater importance to maintaining simulation stability. As a consequence, the agent tends to apply filtering more consistently, adjusting the filter strength dynamically rather than switching it off. The inclusion of physics-based terms in the reward systematically promotes stabilization, preventing the growth of numerical instabilities.

\begin{figure}[htpb!]
    \centering
    \includegraphics[width=\linewidth]{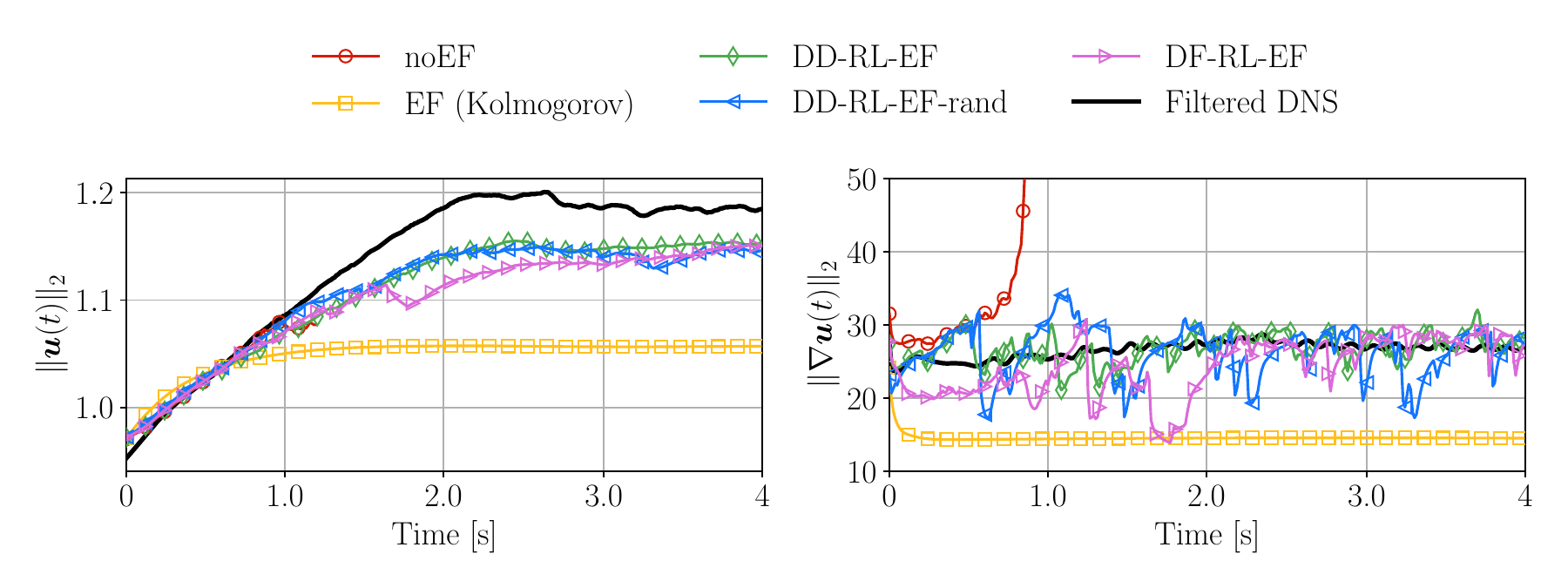}
    \caption{\emph{Test case 1}. Time evolution of the velocity and velocity gradients' norms, for the filtered DNS, noEF, EF (Kolmogorov), and for the novel RL-EF methods.}
    \label{fig:cylinder-norms}
\end{figure}

\begin{figure}[htpb!]
    \centering
    \includegraphics[width=\linewidth]{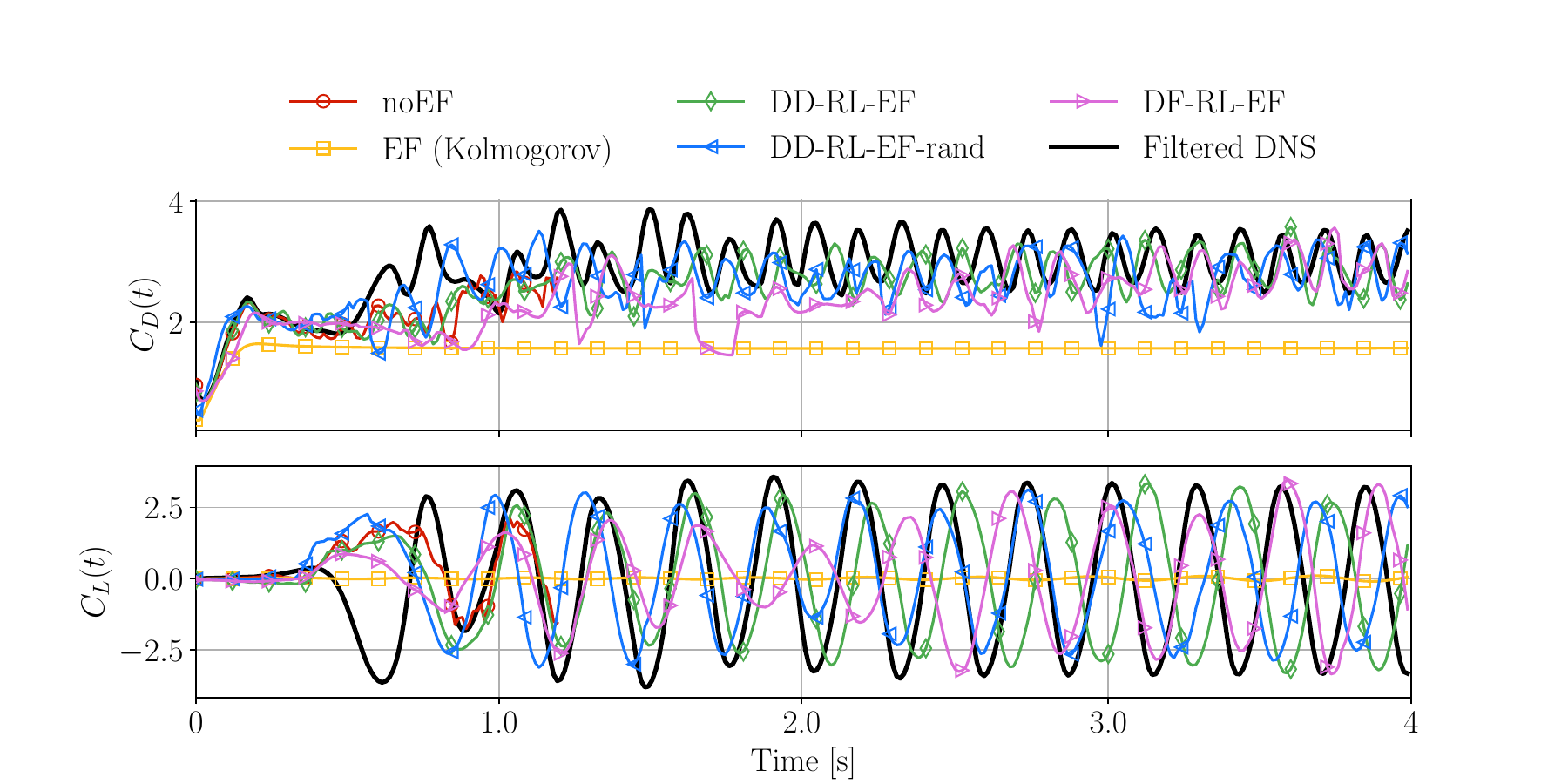}
    \caption{\emph{Test case 1}. Time evolution of drag and lift coefficients, for the filtered DNS, noEF, EF with $\delta=\eta$, and for the novel RL-EF methods. Drag and lift forces are computed on the cylinder boundary.}
    \label{fig:cylinder-forces}
\end{figure}

The above discussion is confirmed by Figures \ref{fig:cylinder-norms} and \ref{fig:cylinder-forces}, which display a comparison between the noEF, RL-EF, EF (Kolmogorov) and the reference filtered DNS. In particular, Figure \ref{fig:cylinder-norms} compares the methods in terms of the velocity norms $\|\bu\|_2$ and $\|\nabla \bu\|_2$, while Figure \ref{fig:cylinder-forces} analyzes the drag and lift forces acting on the cylinder. The definition of the lift and drag adimensional coefficients is:
\begin{equation}
\begin{split}
    &C_D(t) = \dfrac{2}{U^2 L}\int_{\partial\Omega_C}((2\nu \nabla \bu(t) -p (t)\mathbf{I}) \cdot \boldsymbol{n}_C) \cdot \boldsymbol{t}_C \, ds, \\
    &C_L(t) = \dfrac{2}{U^2 L}\int_{\partial\Omega_C}((2\nu \nabla \bu(t) -p(t)\mathbf{I}) \cdot \boldsymbol{n}_C) \cdot \boldsymbol{n}_C \, ds,
\end{split}
\end{equation}
where $U$ and $L$ are the characteristic velocity and length of the test case, while $\partial \Omega_C$ is the cylinder boundary, $\boldsymbol{n}_C$ and $\boldsymbol{t}_C$ are the unitary vectors in the normal and tangential direction to the cylinder.

Figures \ref{fig:cylinder-norms} and \ref{fig:cylinder-forces} show that the RL-EF methods are more accurate and closer to the filtered DNS reference both in $\|\bu\|_2$ and in gradients $\|\nabla \bu\|_2$, than noEF and standard EF.
Across all the considered metrics, the over-diffusive nature of the standard EF (Kolmogorov) approach is evident in the absence of oscillations and in an almost flat temporal behaviour. While this guarantees stability, it also suppresses relevant dynamical features of the flow.
The noEF simulation, on the other hand, becomes unstable and blows up after $t=1$. For this reason, its corresponding red curve is not reported beyond this time in the plots.
All the RL-EF strategies successfully recover the correct energy evolution and closely approach the reference solution. However, the velocity gradients exhibit temporal oscillations around the reference values. This behaviour is not unexpected and is intrinsic to this class of simulations. Notably, similar oscillations were also observed in~\cite{Ivagnes2025}, even though in that work each action $\delta_n$ was computed as the optimal value obtained from a separate optimization performed independently at every time step $t_n$. Therefore, the presence of oscillations cannot be attributed to suboptimality of the learned policy. In fact, in~\cite{Ivagnes2025} the filter parameter was determined using full access to the reference data and via a time-local optimization at each step, which represents a significantly stronger setting. In the present work, instead, the RL agent learns a single policy from limited training data and is then required to generalize in a time-extrapolation regime. Observing similar oscillatory patterns under these more challenging conditions further supports the robustness and effectiveness of the proposed RL framework.
Additionally, the data-driven methods closely approach the reference lift and drag coefficients, both in the amplitude and in the frequency. The data-free strategy improves EF (Kolmogorov), but for some time steps in $[1.5, 2.5]$ slightly underestimates the forces amplitudes and the velocity norms.

\begin{figure}[htpb!]
    \centering
    \includegraphics[width=\linewidth]{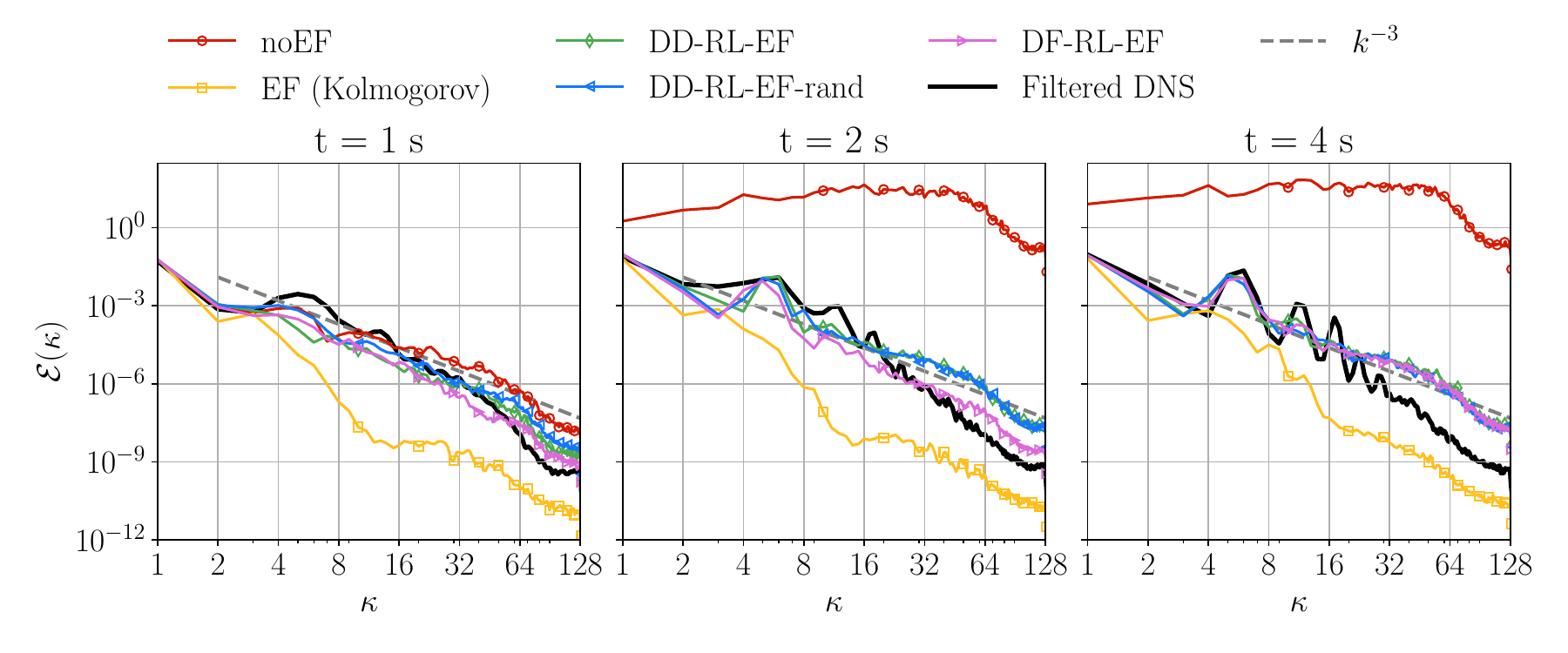}
    \caption{\emph{Test case 1}. Energy spectrum at three time instances, for the filtered DNS, noEF, EF with $\delta=\eta$, and for the novel RL-EF methods. The spectrum is computed at the cylinder wake.}
    \label{fig:cylinder-spectra}
\end{figure}

Figure \ref{fig:cylinder-spectra} displays the energy spectrum at different wavenumbers $\mathcal{E}(\kappa)$ for the proposed methods, compared with the theoretical turbulence trend $\kappa$ in two-dimensional test cases.
At $t=1$, when the vortex shedding is not yet fully developed, all methods exhibit a comparable large-scale energy distribution. However, EF (Kolmogorov) approach already shows a pronounced damping of the high-wavenumber modes, resulting in a premature energy decay in the small scales. This confirms its over-diffusive character.
At later times ($t=2$ and $t=4$), the differences become more evident. The noEF simulation accumulates excessive energy at intermediate and high wavenumbers, as it blows up. In contrast, EF (Kolmogorov) strongly suppresses small-scale activity, producing spectra that are overly steep and significantly below the reference solution in the dissipative range.
All RL-EF strategies achieve a substantially improved balance. They correctly preserve the large-scale energy content while controlling the high-wavenumber range without introducing excessive diffusion. The resulting spectra remain close to the reference slope, including in the inertial range where the $\kappa^{-3}$ behaviour is expected. 
Moreover, at $t=1$ and $t=2$, the DF-RL-EF spectrum exhibits lower energy at large wavenumbers (i.e., small scales) and remains closer to the filtered DNS compared to the DD-RL-EF strategies. This behaviour is directly consistent with the larger $\delta_n$ values observed in Figure \ref{fig:cylinder-actions}, indicating that the data-free agent applies stronger filtering during these stages.
This confirms that the physics-based reward effectively promotes stabilization by penalizing the uncontrolled growth of small-scale energy. In particular, it enables the agent to regulate the dissipative range in a physically meaningful way, preventing energy pile-up at high wavenumbers while preserving the correct large-scale dynamics. The result is a more controlled small-scale behaviour without introducing the excessive diffusion characteristic of the standard EF (Kolmogorov) approach.
Overall, this demonstrates that a purely physics-informed reward is sufficient to guide the agent toward a spectrally consistent filtering strategy, capable of balancing stability and accuracy across scales.

\begin{figure}[htpb!]
    \centering
    \subfloat{\includegraphics[width=.85\linewidth, trim={15cm 0 0 0}, clip]{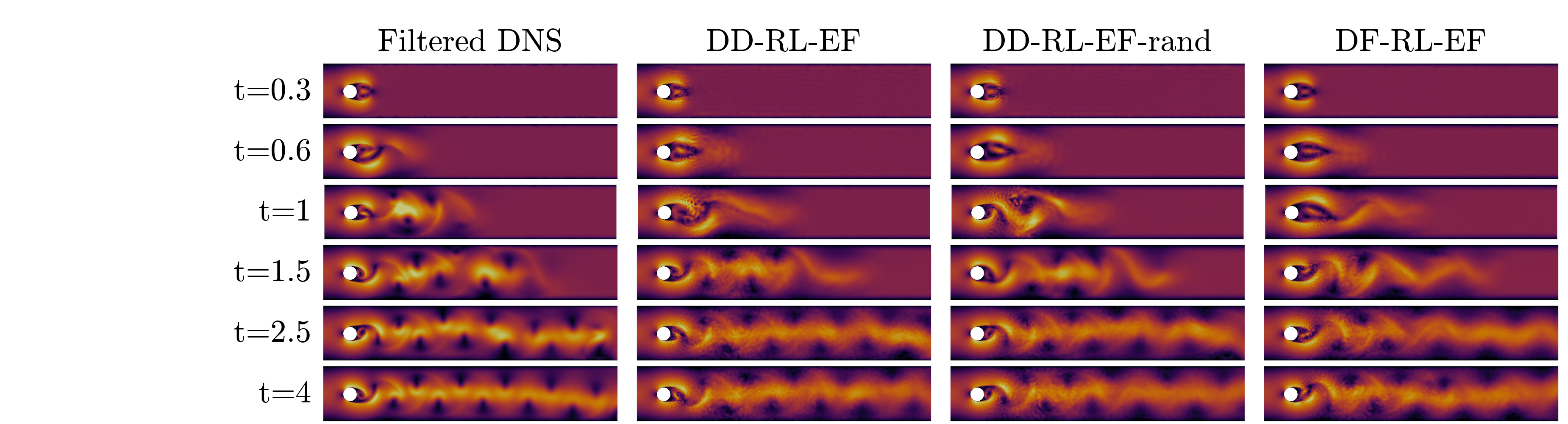}}
    \subfloat{\includegraphics[width=.15\linewidth, trim={90cm 28cm 25cm 10cm}, clip]{images_cylinder/legend.png}}
    \caption{\emph{Test case 1}. Velocity magnitude fields for the RL-EF methods, and for the filtered DNS reference. The fields are represented at different time instances.}
    \label{fig:cylinder-RL-graphics}
\end{figure}

Furthermore, Figure \ref{fig:cylinder-RL-graphics} shows the velocity magnitude fields across different time instances. Compared to noEF and EF (Kolmogorov) in Figure \ref{fig:cylinder-pre-graphics}, all RL-EF are qualitatively able to accurately retrieve the vortex shedding, without being either overdiffusive or noisy.

\begin{figure}[htpb!]
    \centering
    \includegraphics[width=0.9\linewidth]{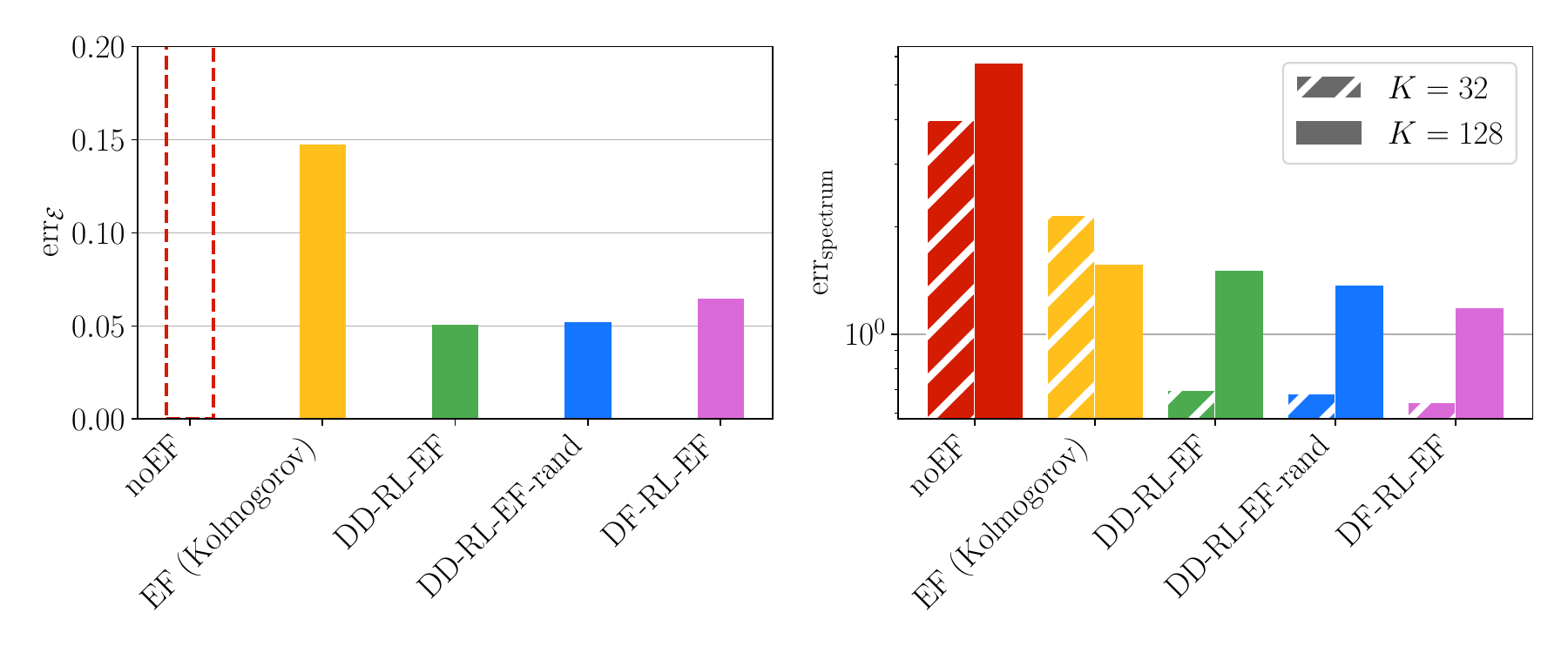}
    \caption{\emph{Test case 1}. Error of the kinetic energy averaged in time ($\mathrm{err}_{\mathcal{E}}$), and log-space error of the energy spectrum averaged with respect to the wavenumbers and in time ($\mathrm{err}_{\text{spectrum}}$). The energy spectrum error is computed over two different wavenumbers ranges, namely $K=32$ and $K=128$.}
    \label{fig:glob-errs-cyl}
\end{figure}

Finally, we compare the methods in terms of global metrics, obtained averaging in time the accuracy of the RL-EF kinetic energy, and of the RL-EF energy spectrum at different scales. 
We consider \textit{global} quantities, e.g., energy and enstrophy, because the classical relative error in the $L^2$ norm is not appropriate to quantify the performance of the proposed methodology. As noticed in~\cite{Ivagnes2025, strazzullo2024variational}, it might indeed happen that EF-like and noEF solutions have similar relative $L^2$ error, even if the results are different.
In particular, we consider:  

\begin{equation}
    \mathrm{err}_{\mathcal{E}}=\dfrac{1}{N}\sum_{n=1}^N \left| \dfrac{\mathcal{E}_n - \mathcal{E}_n^{\text{ref}}}{\mathcal{E}_n^{\text{ref}}} \right|, \quad 
    \mathrm{err}_{\text{spectrum}}=\dfrac{1}{N}\sum_{n=1}^N \left(\dfrac{1}{K} \sum_{\kappa=1}^{K} \log_{10}{\left( \dfrac{\mathcal{E}_n(\kappa)}{\mathcal{E}_n^{\text{ref}}(\kappa)} \right)} \right),
    \label{eq:err-global}
\end{equation}
where $K$ is the maximum wavenumber considered to compute the error. In $\mathrm{err}_{\text{spectrum}}$, we consider a $\log_{10}$-scaled error, as it better describes the relative magnitude of discrepancies across several orders of magnitude and emphasizes differences in the low-energy, large wavenumber range of the spectrum, which are often crucial for assessing the fidelity of numerical schemes in turbulent and multiscale flows.
Figure \ref{fig:glob-errs-cyl} displays the above-mentioned metrics for all the RL-EF solutions, in comparison with noEF and EF (Kolmogorov). Moreover, we consider $\mathrm{err}_{\text{spectrum}}$ for two different $K$ values: the biggest one, retaining all high-energy (small wavenumbers) and low-energy (large wavenumbers) scales, and a smaller one retaining only the large-energy scales. 
This allows us to identify the accuracy of the RL-EF both on the dominant energy-containing scales (small $K$), and on the entire spectral range (large $K$).

For the cylinder test case (Figure \ref{fig:glob-errs-cyl}), regarding the energy error $\mathrm{err}_{\mathcal{E}}$, all RL-EF methods are considerably more accurate ($5\%$ error) then EF (Kolmogorov) ($15\%$ error). The red line for noEF simulation is dashed as its solution blows up. 
Regarding the energy spectrum error, the RL-EF methods, and in particular the reference-free approach DF-RL-EF, are the most accurate for both $K$ values. However, there is a larger difference in the error when only the large scales are considered ($K=32$). This indicates that the choice of the filtering strategy has a stronger impact on the accuracy of the energy distribution in the dominant structures of the flow.
In particular, EF (Kolmogorov) exhibits a larger error in the large-scale energy content, as it is not able to recover the physical vortex shedding. In contrast, the RL-based strategies recover the relevant flow dynamics and better preserve the energy distribution across scales. When the full spectrum is considered ($K=128$), the differences between the RL approaches remain limited, with DF-RL-EF achieving the lowest overall error. EF (Kolmogorov) has a similar performance to the RL-EF approaches because it does not exhibit spurious oscillations or noisy behaviour at the small scales.
In conclusion, it is important to emphasize that:
\begin{itemize}
    \item In the \textbf{reference-guided} approaches, reference data are only provided within a limited time window, namely $[0,1]$ (and $[0,1.4]$ for DD-RL-EF-rand). In this interval, vortex shedding is not yet fully developed. Despite this restricted and pre-transient training regime, the RL agent successfully learns how to modulate the filtering intensity from the velocity state alone. Remarkably, the learned policy generalizes to the time-extrapolation regime, where it maintains both stability and accuracy of the RL-EF simulations even after the flow becomes fully developed.
    \item In the \textbf{reference-free} setting, the RL agent is trained without any access to reference data and relies solely on physics-based information. Nevertheless, it autonomously learns an effective filtering strategy that preserves stability and accuracy. The resulting performance is very close to that of the data-driven approaches, demonstrating that physically informed rewards alone are sufficient to guide the learning process toward robust and high-quality solutions, drastically reducing the computational time needed for the training, since no reference is used and no DNS is needed, as deeper analyzed in Section \ref{comp-details}.
\end{itemize}
Overall, these results highlight the robustness and generalization capability of the proposed RL framework, showing that effective stabilization policies can be learned even from limited or entirely absent reference data.

\subsection{Test case 2: decaying homogeneous turbulence}
\label{sec:test-case-2}
\begin{figure}[htpb!]
    \centering
    \subfloat[Coarse mesh]{\includegraphics[width=0.2\linewidth, trim={15cm 1cm 15cm 0cm}, clip]{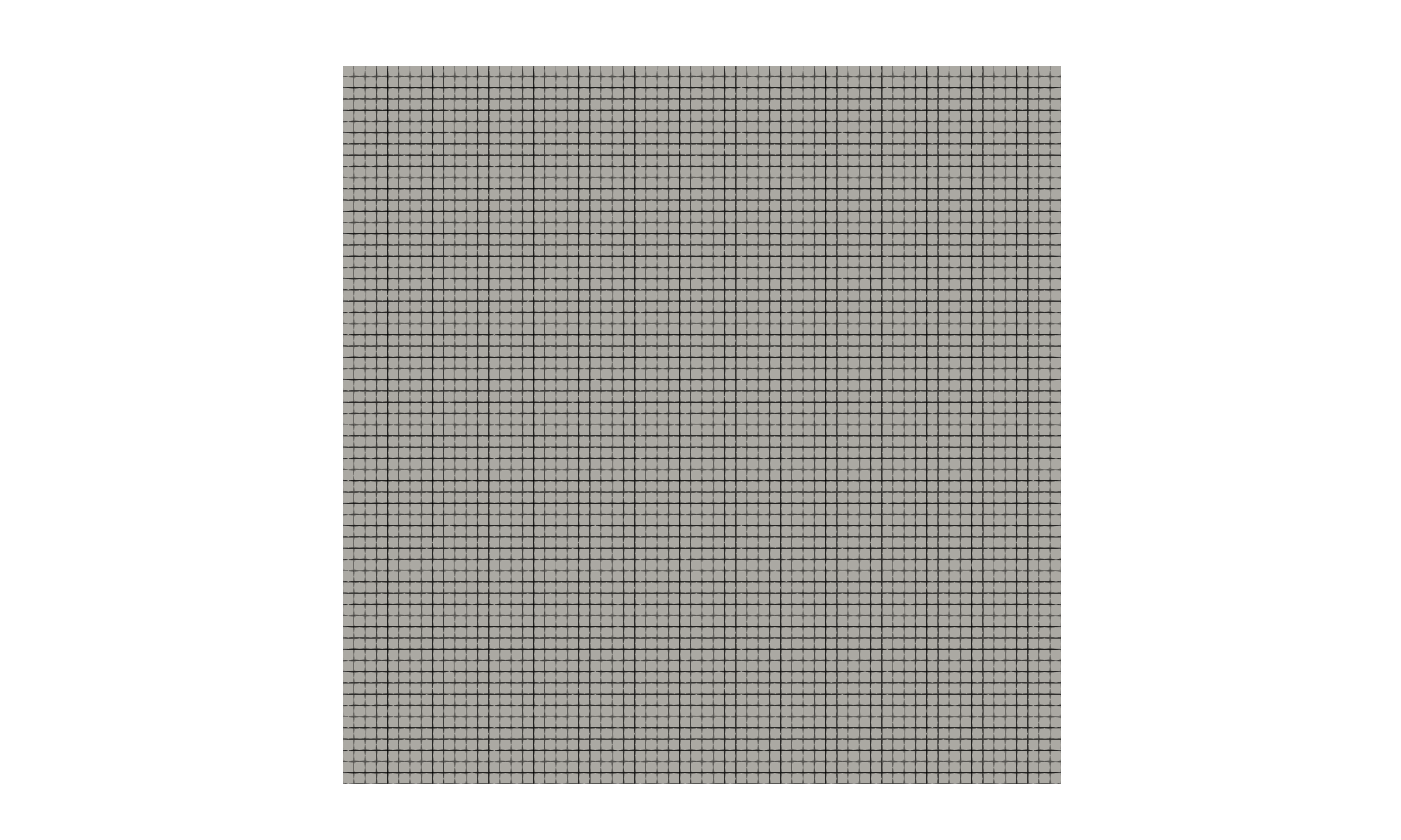}
    \label{coarse-decaying}}
    \subfloat[Fine mesh]{\includegraphics[width=0.2\linewidth, trim={15cm 1cm 15cm 0cm}, clip]{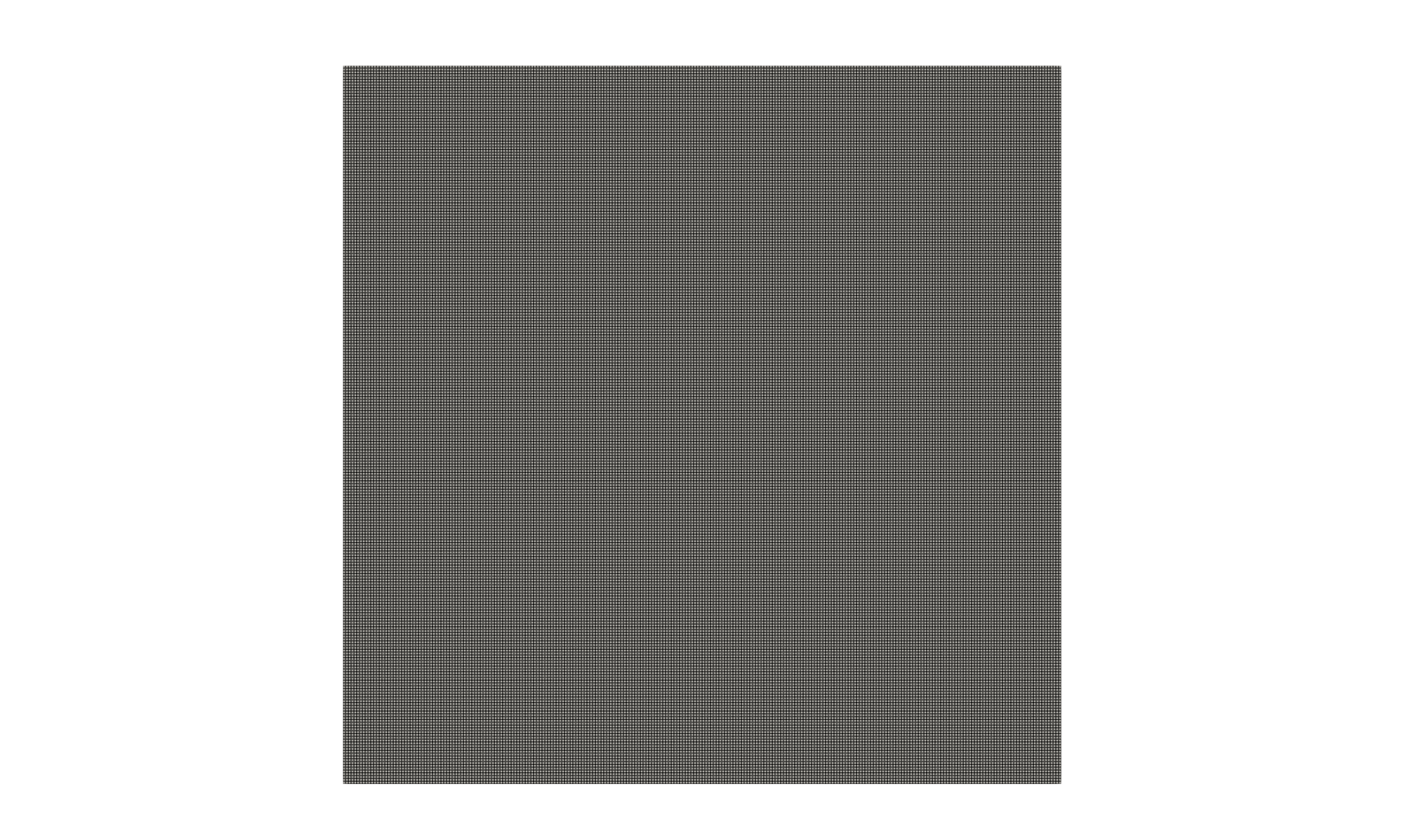}\label{fine-decaying}}
    \subfloat[Velocity magnitude at $t=0$][\\Velocity magnitude at $t=0$]{\includegraphics[width=0.27\linewidth, trim={5cm 5cm 25cm 5cm}, clip]{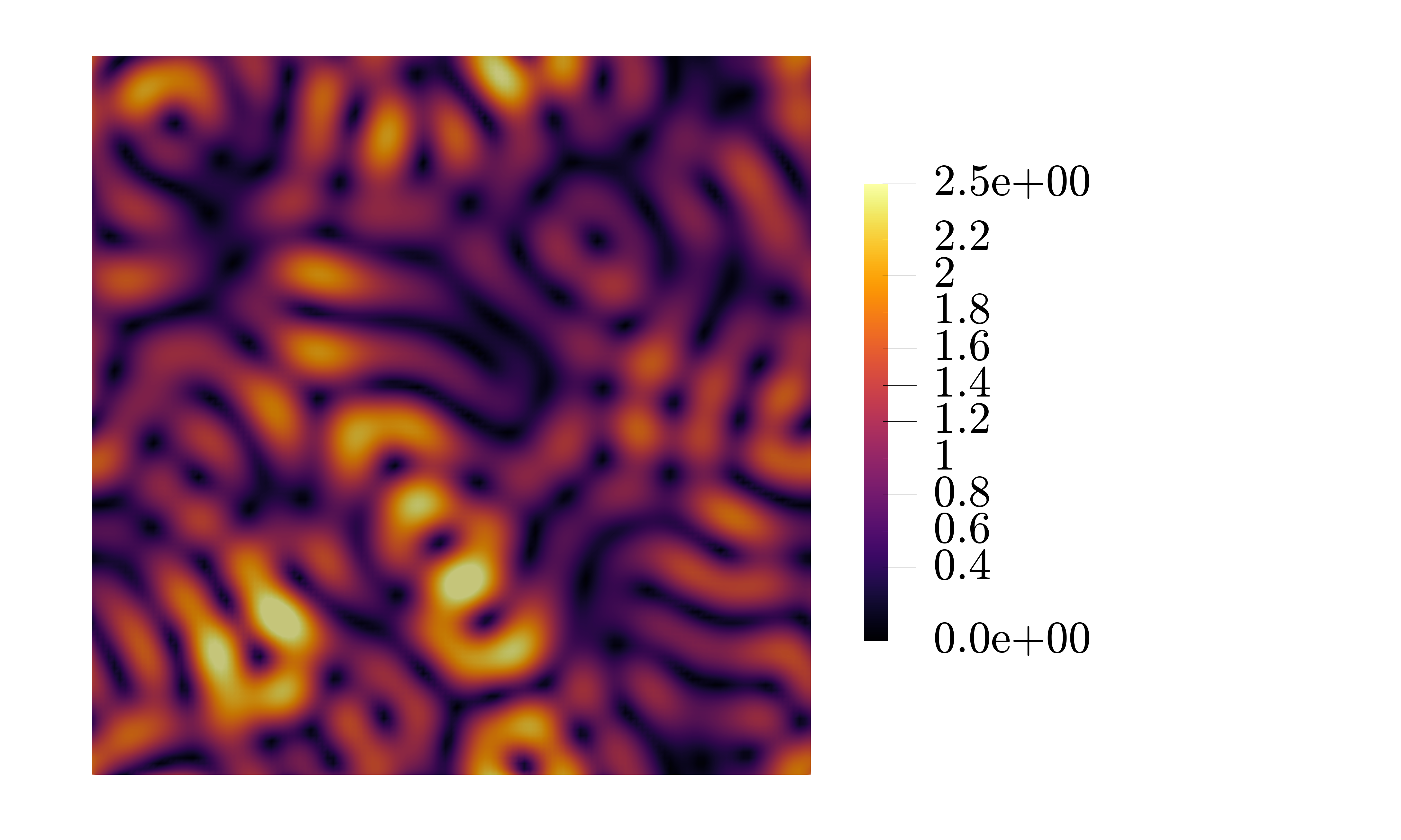}\label{init-vel-decaying}}
    \subfloat[Vorticity at $t=0$][\\Vorticity at $t=0$]{\includegraphics[width=0.27\linewidth, trim={5cm 5cm 25cm 5cm}, clip]{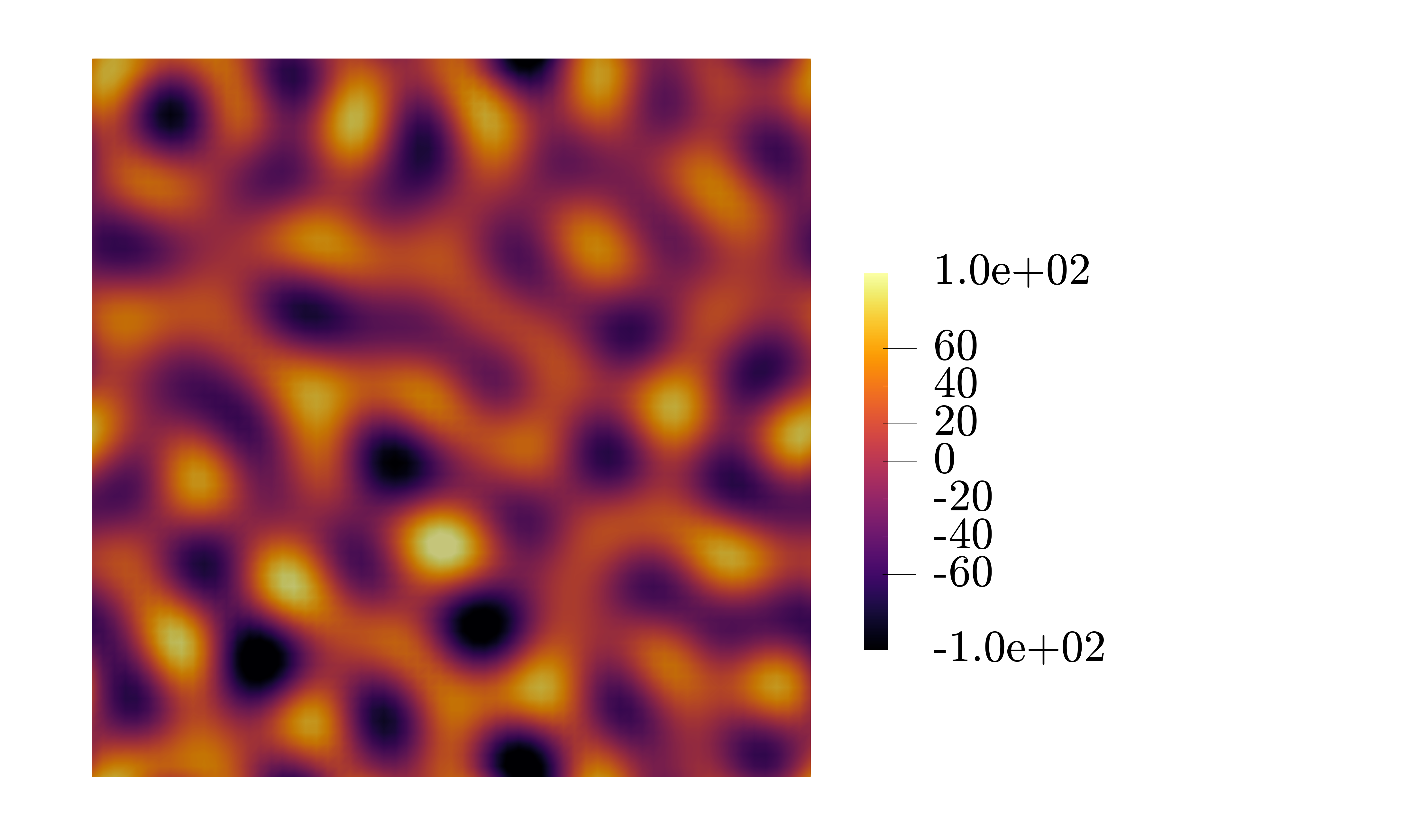}\label{init-vort-decaying}}
    \caption{\emph{Test Case 2}. The two grids used for the decaying homogeneous turbulence test case \protect \subref{coarse-decaying} and \protect \subref{fine-decaying}, the initial velocity \protect \subref{init-vel-decaying} and vorticity \protect \subref{init-vort-decaying} fields. The coarse grid is used for noEF, EF, and RL-EF simulations, while the fine mesh is used for the DNS reference simulation.}
    \label{fig:meshes-decaying}
\end{figure}

The second test case is the two-dimensional decaying homogeneous turbulent flow in a square at $Re=\num{40000}$.
Table \ref{tab:dofs-decaying} reports the number of degrees of freedom used for the velocity ($N_h^{\boldsymbol{u}}$) and for the pressure ($N_h^p$) field, the grid resolution for the coarse and refined grids, and the corresponding mesh sizes. Both grids are equally-spaced as represented in Figure \ref{coarse-decaying} and \ref{fine-decaying}, with periodic boundary conditions and quadrilateral elements.

\begin{table}[H]
    \caption{\emph{Test Case 2}. Number of degrees of freedom, minimum and maximum sizes for the grids considered.}
    \centering
    \begin{tabular}{ccccc}
    \toprule
    & $N_h^{\boldsymbol{u}}$ & $N_h^p$ & Grid resolution& $h_{min}=h_{max}$\\
    \midrule
        Coarse mesh & $16641$ & $4225$& $64\times 64$&$1.56e-2$  \\
         Fine mesh & $263169$ & $66049$ &$256 \times 256$&$3.91e-3$ \\
    \bottomrule     
    \end{tabular}
    \label{tab:dofs-decaying}
\end{table}

As in \cite{Agdestein2025MLLES}, the velocity field is initialized with the following procedure:

\begin{enumerate}
    \item The velocity is sampled in spectral space to match a prescribed energy $\hat{\mathcal{E}}_{\kappa}$ at each wave number $\kappa$. Each component of the spectral velocity is indeed initialized such that $\|\hat{\boldsymbol{u}}_{\kappa} \|= |\mathrm{a}_{\kappa}|$, where $\mathrm{a}_{\kappa}=\sqrt{2\hat{\mathcal{E}}_{\kappa}} \text{e}^{2\pi i \tau_{\kappa}}$ and $\tau_{\kappa}$ is a random phase shift. 
    \item The velocity is then projected to make it divergence-free: $\hat{\boldsymbol{u}}_{\kappa}=\frac{\mathrm{a}_{\kappa} \hat{\boldsymbol{P}}_{\kappa} \boldsymbol{e}_{\kappa}}{\hat{\boldsymbol{P}}_{\kappa} \boldsymbol{e}_{\kappa}}$, where $\hat{\boldsymbol{P}}_{\kappa}$ is a projector: $\hat{\boldsymbol{P}}_{\kappa}=I-\frac{\kappa \kappa^T}{\kappa^T \kappa}$ and $\boldsymbol{e}_{\kappa}$ is a random unit vector. Since the term $\frac{\kappa \kappa^T}{\kappa^T \kappa} \boldsymbol{e}_{\kappa}$ provides the projection of $\boldsymbol{e}_{\kappa}$ onto $\kappa$, applying the operator $\hat{\boldsymbol{P}}_{\kappa}$ leads to $\kappa^T (\hat{\boldsymbol{P}}_{\kappa} \boldsymbol{e}_{\kappa})=0$, and hence $\kappa^T \hat{\boldsymbol{u}}_{\kappa}=0$. This condition is equivalent to the incompressibility constraint in the frequency domain. In 2D $\boldsymbol{e}_{\kappa}=(\cos{(\theta_{\kappa})}, \sin{(\theta_{\kappa})})$ with $\theta_{\kappa} \sim \mathcal{U}[0, 2\pi]$.
    \item The velocity field is finally obtained by a transformation into physical space ($\boldsymbol{u}_1 = \text{FFT}^{-1}(\hat{\boldsymbol{u}})$), and a further projection step such that $\boldsymbol{M}_h \boldsymbol{u}_1=\boldsymbol{0}$, where $\boldsymbol{M}_h$ is the discrete divergence operator.
\end{enumerate} 

\begin{figure}[htpb!]
    \centering \subfloat{\includegraphics[width=.8\linewidth]{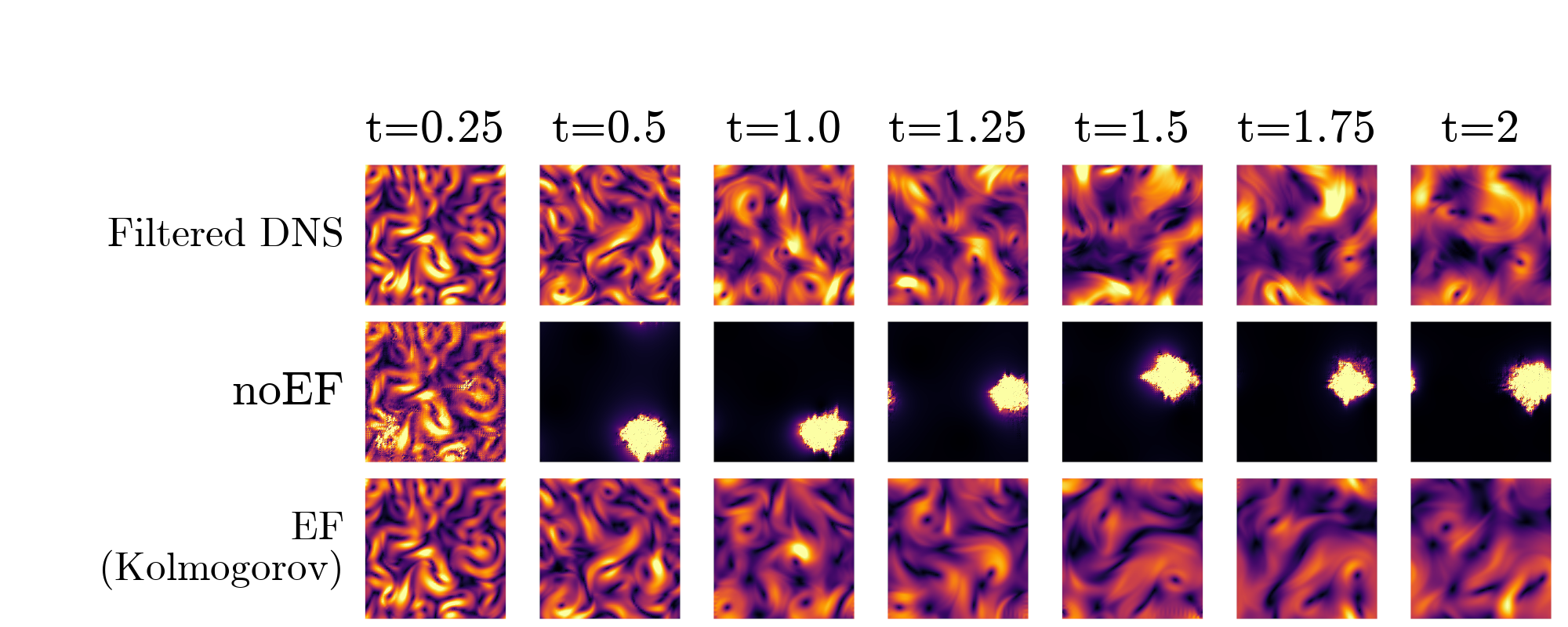}}
    \subfloat{\includegraphics[width=.2\linewidth,
    trim={110cm 20cm 5cm 5cm}, clip]{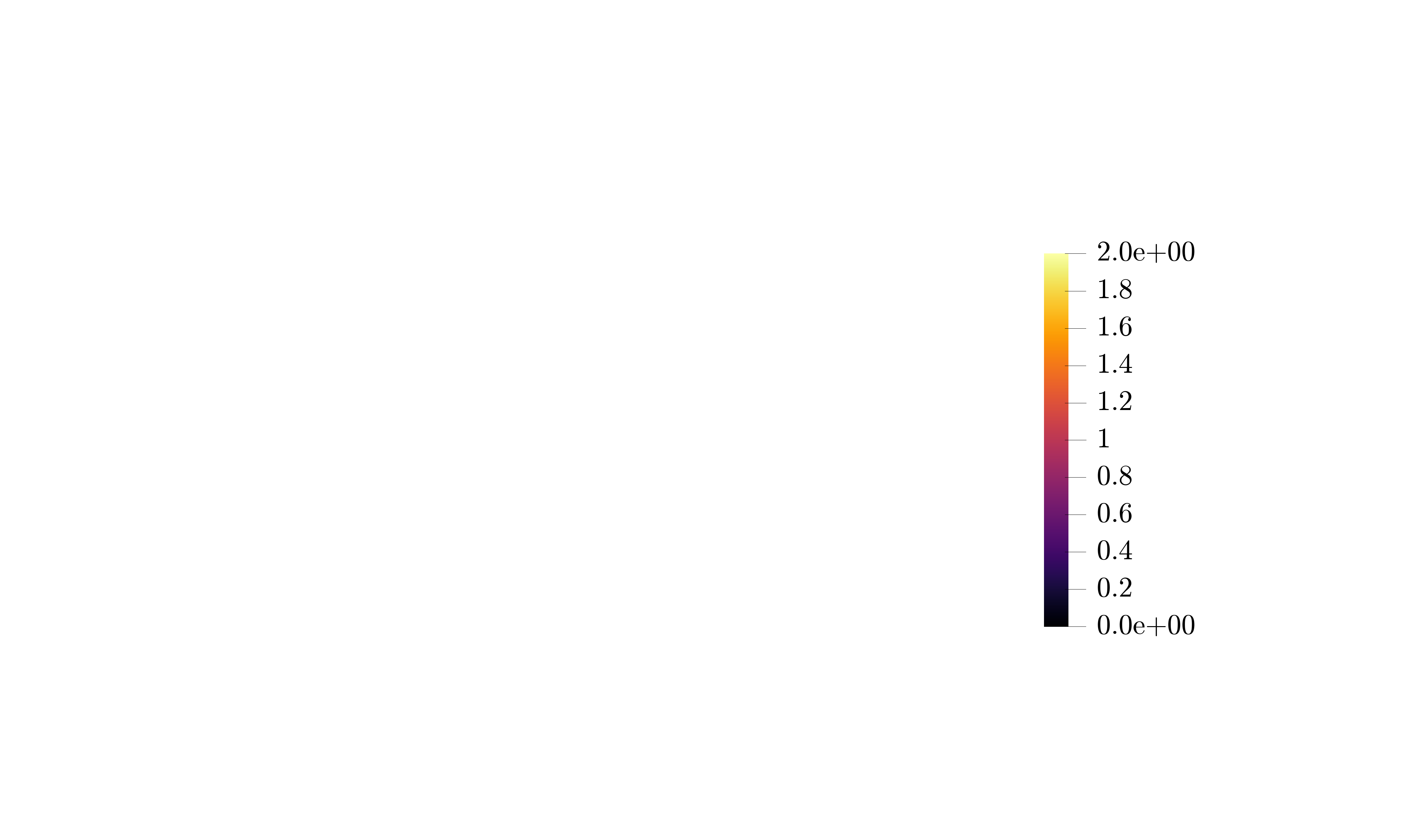}}
    \caption{\emph{Test Case 2}. Velocity magnitude fields for the filtered DNS, the noEF and EF with $\delta=\eta$ simulations, for the decaying turbulence test case. The fields are represented at different time instances.}
    \label{fig:decaying-pre-graphics}
\end{figure}

Each simulation lasts $T=2$ seconds, with a time step $\Delta t=\num{1e-3}$, for a global number of time steps $N=2000$.
As a preliminary assessment, Figure~\ref{fig:decaying-pre-graphics} shows some representative solutions in time of the filtered DNS, noEF, and EF (Kolmogorov) velocity fields.
Similarly to the flow past a cylinder case, the noEF solution becomes unstable and blows up after only a few time steps. In contrast, the EF (Kolmogorov) solution closely follows the filtered DNS up to approximately $t=0.5$. Beyond this time, however, it progressively diverges from the reference solution due to excessive diffusion. This over-diffusive behaviour is evident in the size of the vortical structures, which appear larger compared to those of the filtered DNS.
The filtered DNS reference flow initially consists of small-scale vortical structures that progressively merge and enlarge over time, leading to a characteristic decay of turbulence and a predominantly dissipative behaviour.
This evolution is fundamentally different from the flow past a cylinder. In the cylinder case, vortex shedding develops only after a transient phase, with the kinetic energy first increasing and then reaching a periodic regime. In contrast, in the decaying turbulence scenario, the kinetic energy monotonically decreases as time advances, and the dynamics are entirely governed by the gradual dissipation of structures across scales.
This marked difference raises a crucial question: \\
\textit{can the RL-based filtering strategy adapt to a regime where no sustained energy production occurs, and where the primary challenge is to preserve the correct decay rate without introducing artificial dissipation?}
The following results address precisely this point.

\begin{figure}[htpb!]
    \centering
    \includegraphics[width=\linewidth]{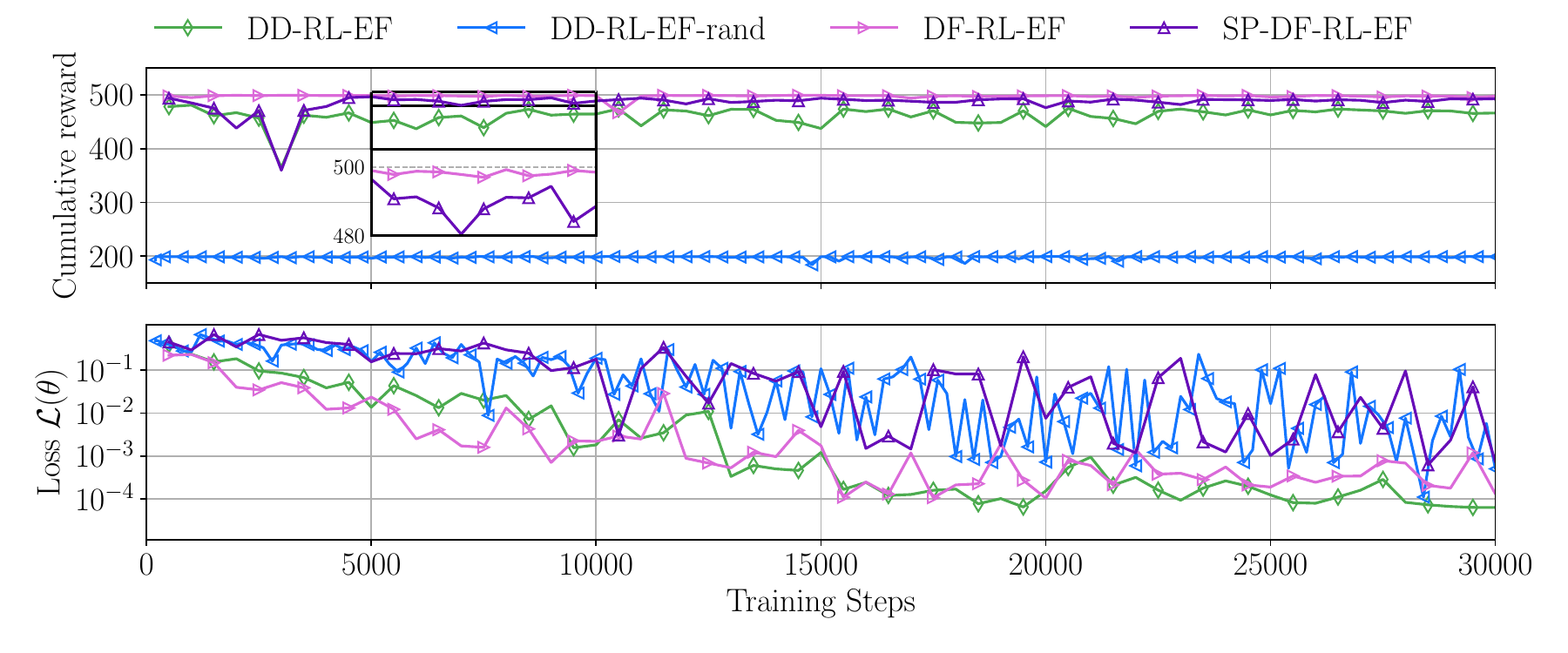}
    \caption{\emph{Test Case 2}. Cumulative reward and loss value during training steps, for the RL-EF methods.}
    \label{fig:decaying-rewards}
\end{figure}

Figure \ref{fig:decaying-rewards} shows the cumulative reward and the DQN loss function across training episodes for different RL-EF strategies.
Being this test case an \textit{energy} and \textit{enstrophy}-preserving test case, we also investigate the \textbf{SP-DD-RL-EF} strategy introduced in Section \ref{sec:rl-ef}.
In this case, the maximum cumulative reward is $N_{\text{train}}=N/4=500$ ($T_{\text{train}}=0.5$) for all RL-EF methods, except for DD-RL-EF-rand, which has shorter episodes of $N_{\text{train}}^{\text{rand}}=N/10=200$ time steps (corresponding to $T_{\text{train}}^{\text{rand}}=0.2$).
All the RL-EF strategies reach a plateau close to their maximum theoretical value, but with different loss decay.
First of all, differently from the first test case, the loss decay is in globally much slower, and it reaches $\sim \num{1e-4}$ only in the cases of DD-RL-EF and DF-RL-EF.
DD-RL-EF-rand and SP-DF-RL-EF exhibit more oscillatory loss behaviour. 
In the DD-RL-EF-rand case, this behaviour is primarily related to the random initialization of each episode, as already observed in the first test case. The variability introduced at the beginning of every training episode increases the dispersion in the learning dynamics.
In the structure-preserving setting, instead, the effect can be attributed to the additional physical constraints embedded in the reward formulation. While these constraints enforce consistency with the underlying physics, they also make the optimization landscape more complex, leading to less stable training and making the identification of the optimal filter parameter more challenging.

\begin{figure}[htpb!]
    \centering
    \includegraphics[width=\linewidth]{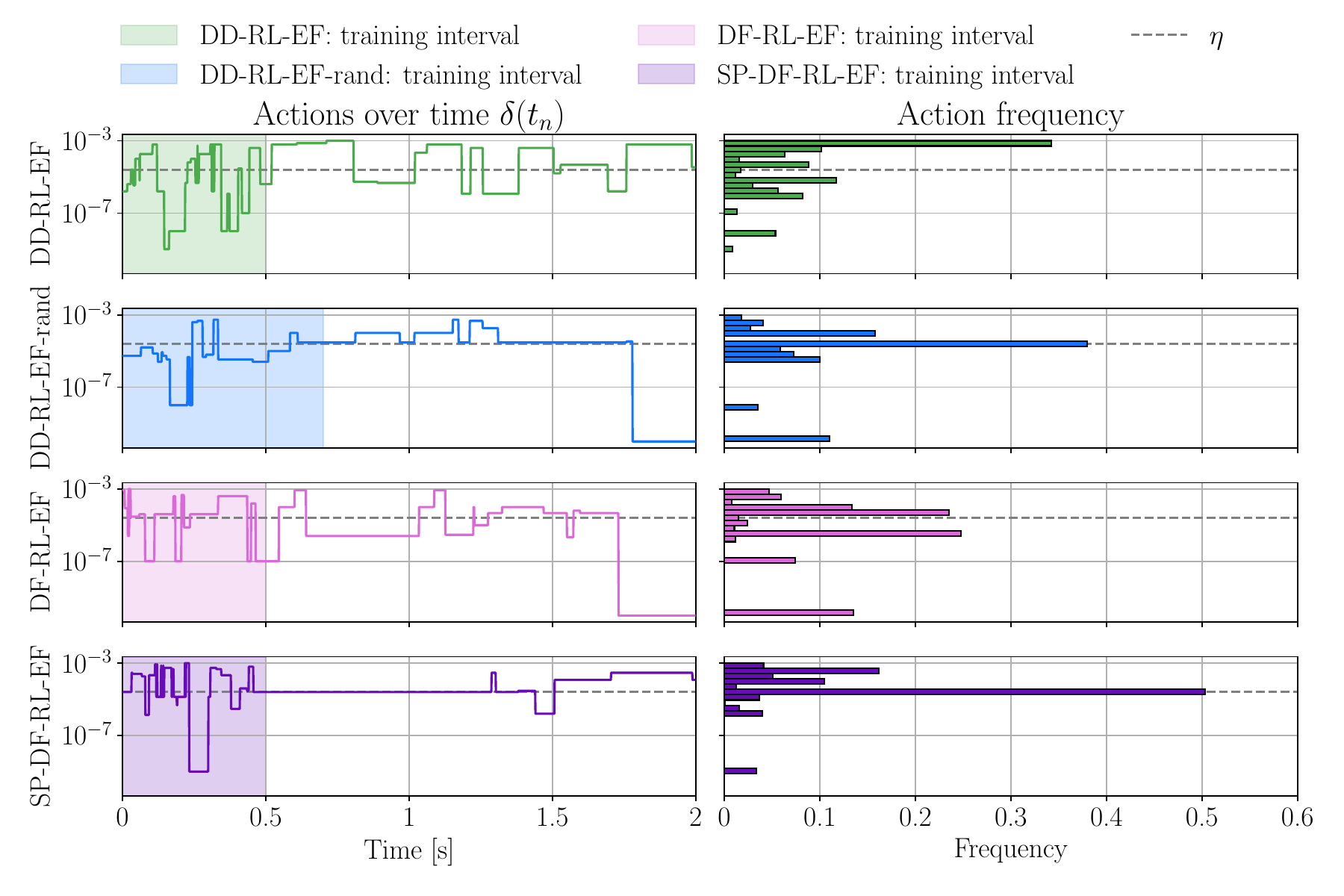}
    \caption{\emph{Test Case 2}. Action values over testing time, and corresponding frequency, for the RL-EF methods.}
    \label{fig:decaying-actions}
\end{figure}

Figure~\ref{fig:decaying-actions} illustrates how the pre-trained RL agents select the filtering actions throughout the simulation, together with the corresponding distribution of the most frequently chosen filter values.
In this test case, the agents exhibit a 
less oscillatory behaviour compared to the previous scenario. In particular, the selected action values are maintained for longer consecutive time intervals.
Furthermore, the action frequency distributions reveal that the most frequently selected filter radius $\delta_n$ at $t_n$ coincides with, or is very close to, the Kolmogorov length scale. This trend is especially pronounced in DD-RL-EF-rand and SP-DD-RL-EF, where the agent selects $\delta_n=\eta$ the $40\%$ and the $50\%$ of the time steps $t_n$, respectively. 
This suggests that, in a purely dissipative regime, the learned policy naturally gravitates toward a physically meaningful small-scale cutoff, balancing stabilization and preservation of the correct decay dynamics.

\begin{figure}[htpb!]
    \centering
    \includegraphics[width=\linewidth]{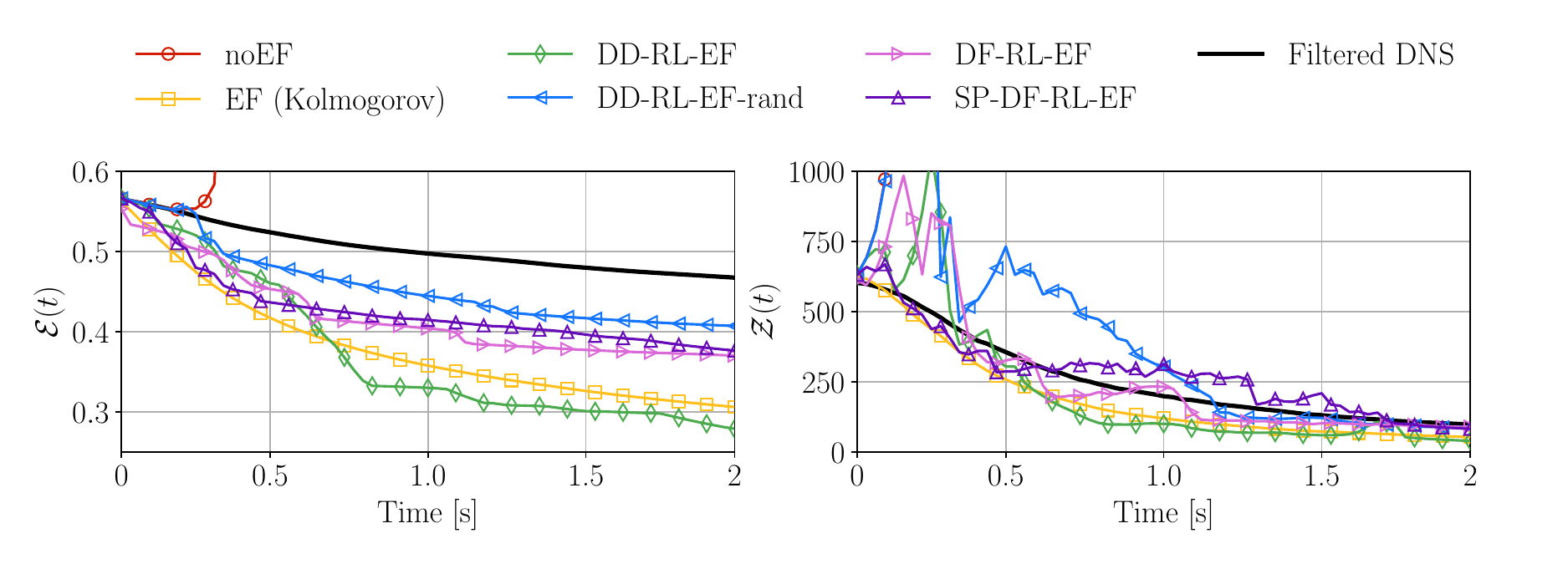}
    \caption{\emph{Test Case 2}. Time evolution of energy and enstrophy, for all the RL-EF simulations, the standard EF, noEF, and the filtered DNS reference simulations.}
    \label{fig:decaying-norms}
\end{figure}

Regarding the time evolution of energy and enstrophy shown in Figure~\ref{fig:decaying-norms}, the noEF simulation becomes unstable and blows up shortly after $t=0.25$. In contrast, the EF (Kolmogorov) approach exhibits a smooth and monotonic decay of both energy and enstrophy, but at the cost of excessive dissipation, leading to an over-damped dynamics compared to the filtered DNS reference.
All RL-EF strategies, with the exception of DD-RL-EF, successfully capture the correct energy decay without introducing the excessive diffusion observed in EF (Kolmogorov). Among them, DD-RL-EF-rand achieves the closest agreement with the filtered DNS in terms of energy evolution. DD-RL-EF, instead, matches the reference well up to $t=0.5$, which corresponds to the training time window, but becomes increasingly over-diffusive at later times. This behaviour is consistent with the larger $\delta_n$ values selected beyond the training horizon, as shown in Figure~\ref{fig:decaying-actions}, indicating limited generalization in the extrapolation regime.
Concerning the enstrophy evolution, DD-RL-EF, DD-RL-EF-rand, and DF-RL-EF exhibit an increase, particularly at early times. This behaviour reflects a temporary amplification of small-scale activity before the dissipative regime becomes dominant. In the SP-DD-RL-EF case, the weak enforcement of enstrophy conservation leads to a more controlled and moderated behaviour. Nevertheless, small enstrophy increases are still observed, which can be attributed to the constraint being imposed in a weak (penalized) form rather than as a hard constraint.

\begin{figure}[htpb!]
    \centering
    \includegraphics[width=\linewidth]{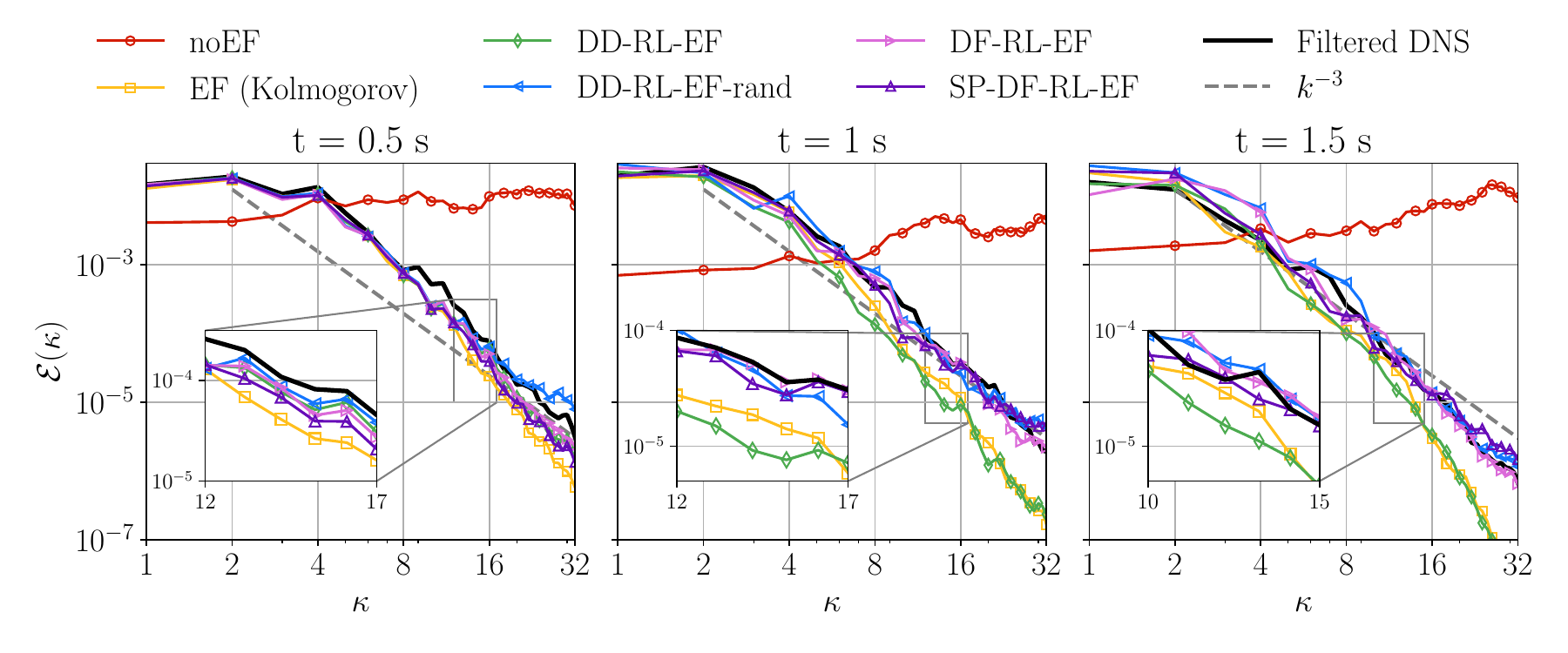}
    \caption{\emph{Test Case 2}. Energy spectrum for the decaying turbulence test case, for all the RL-EF simulations, the standard EF, noEF, and the filtered DNS reference simulations. The spectrum is represented at different time instances.}
    \label{fig:decaying-spectra}
\end{figure}

Figure~\ref{fig:decaying-spectra} compares the energy spectra at $t=0.5, 1$, and $1.5$.
Since the noEF simulation becomes unstable at early times, its spectrum rapidly becomes unphysical, with an excessive accumulation of energy at the high wavenumbers. In contrast, EF (Kolmogorov) systematically underestimates the energy content in the small scales at all reported times, confirming its over-dissipative character.
At $t=0.5$, all RL-EF strategies produce spectra that are very close to each other and in good agreement with the filtered DNS reference. However, at intermediate and large wavenumbers, DD-RL-EF-rand exhibits slightly higher energy levels, indicating residual small-scale oscillations. This behaviour is consistent with the larger enstrophy values observed in Figure~\ref{fig:decaying-norms}, reflecting increased small-scale activity.
At later times ($t \in \{1, 1.5\}$), DD-RL-EF increasingly underestimates the energy in the high-wavenumber range, in line with the stronger dissipation already observed in the energy and enstrophy trends. In contrast, DD-RL-EF-rand, DF-RL-EF, and SP-DD-RL-EF remain in close agreement with the filtered DNS spectrum, accurately reproducing both the large-scale energy content and the small-scale decay.

\begin{figure}[htpb!]
    \centering
    \subfloat{\includegraphics[width=.8\linewidth]{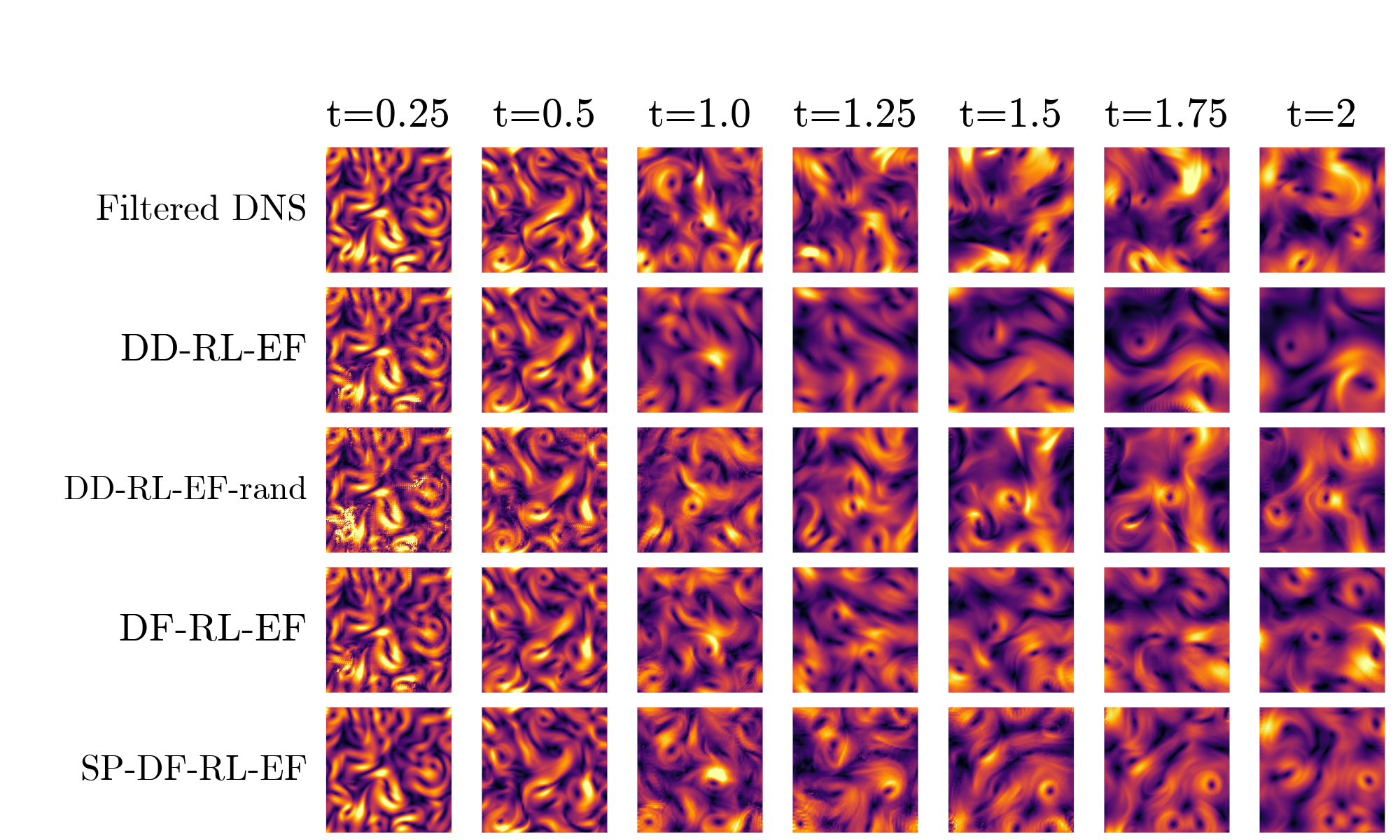}}
    \subfloat{\includegraphics[width=.2\linewidth,
    trim={110cm 5cm 5cm 0cm}, clip]{images_decayingturbulence/legend_decaying.png}}
    \caption{\emph{Test Case 2}. Velocity magnitude fields for the filtered DNS reference and all the RL-EF simulations, for the decaying turbulence test case. The fields are represented at different time instances.}
    \label{fig:decaying-rl-graphics}
\end{figure}

From a qualitative perspective, Figure~\ref{fig:decaying-rl-graphics} compares the RL-EF velocity fields with the filtered DNS at different times. The visual inspection confirms the trends discussed above.
At $t>0.5$, DD-RL-EF clearly exhibits overly smoothed vortical structures, consistent with its over-dissipative behaviour observed in the energy decay and spectral analysis. Conversely, DD-RL-EF-rand shows a noisier small-scale structure at early times (e.g., $t=0.25$), in agreement with its higher enstrophy levels and slightly increased high-wavenumber energy.
The data-free approaches, instead, provide a more balanced solution. The vortical structures remain well-resolved without exhibiting spurious oscillations or excessive smoothing, reflecting a more effective compromise between numerical stabilization and preservation of the physical decay dynamics.

\begin{figure}[htpb!]
    \centering
    \includegraphics[width=0.9\linewidth]{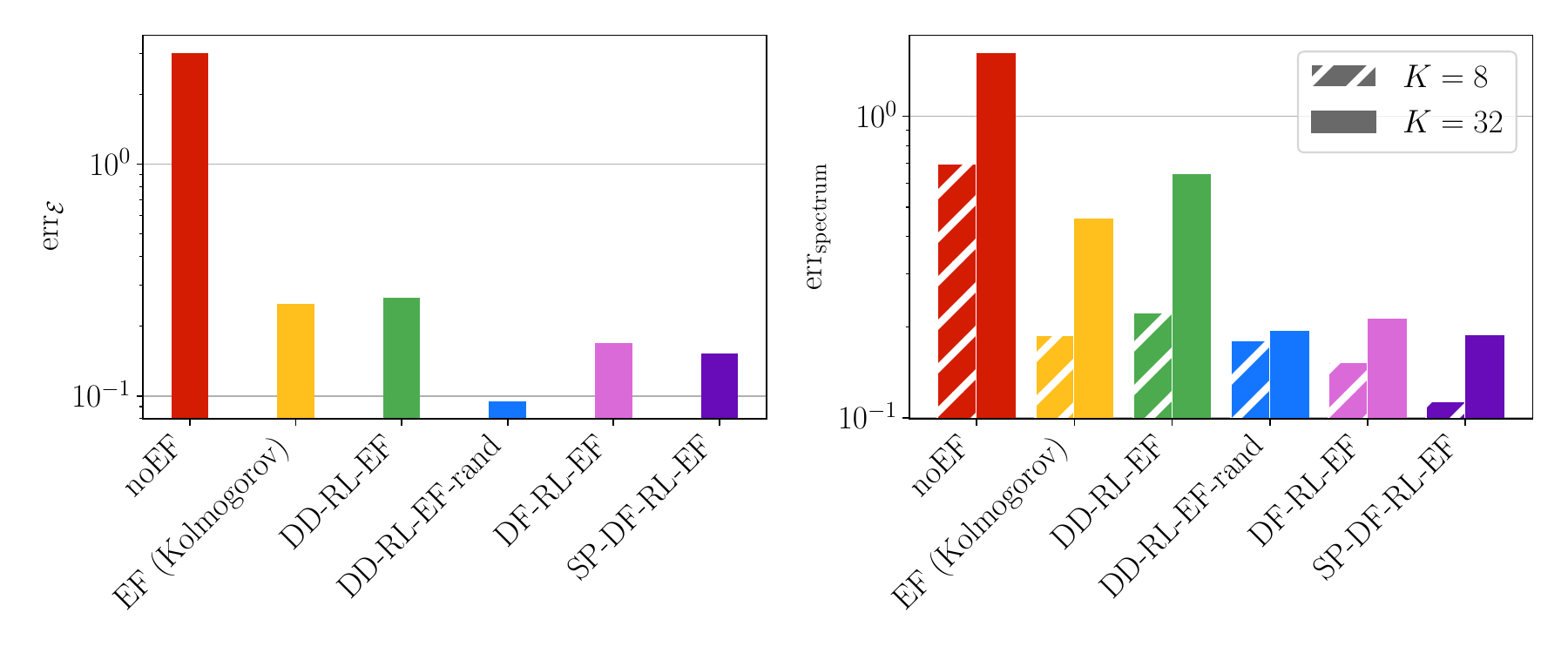}
    \caption{\emph{Test Case 2}. Error of the kinetic energy averaged in time ($\mathrm{err}_{\mathcal{E}}$), and log-space error of the energy spectrum averaged with respect to the wavenumbers and in time ($\mathrm{err}_{\text{spectrum}}$). The energy spectrum error is computed over two different wavenumbers ranges, namely $K=8$ and $K=32$.}
    \label{fig:glob-errs-decaying}
\end{figure}

As in Test Case 1, we conclude the analysis comparing the methods in terms of the global errors on the kinetic energy and on the energy spectrum (the expressions in Equations \eqref{eq:err-global}).
In terms of the error $\mathrm{err}_{\mathcal{E}}$, Figure \ref{fig:glob-errs-decaying} confirms that DD-RL-EF-rand is the most accurate approach, followed by the reference-free methods (SP-DF-RL-EF and DF-RL-EF). As DD-RL-EF is overdiffusive, its energy error is close to the EF (Kolmogorov) one.
Regarding the error in the energy spectrum, we distinguish also in this case the effect of the large scales (maximum wavenumber $K=8$) and of the entire range of scales ($K=32$).
Unlike Test Case 1, here, the discrepancy between the methods becomes more pronounced when the full range of scales is considered, whereas the differences are less significant when only the largest scales are included.
This behaviour reflects the intrinsic difficulty of this test case, where the main challenge lies in accurately controlling the velocity gradients, while avoiding overdissipation.
When considering the full range of scales, the results show a marked difference between the overdiffusive approaches, namely EF (Kolmogorov) and DD-RL-EF, and all the others (DD-RL-EF-rand, DF-RL-EF, and SP-DF-RL-EF). Overdiffusion leads indeed to smaller energy at the large wavenumbers, leading to larger errors at $K=32$.
When a smaller $K$ is considered, all the RL-EF approaches have similar error to EF (Kolmogorov), except for the structure-preserving strategy (SP-DF-RL-EF). The weak constraint in the reward allows to recover the energy spectrum at the dominant scales. This leads to a better agreement with the filtered DNS spectrum, even without training the algorithm with the reference DNS data.

Overall, this test case demonstrates that the proposed RL-based filtering framework is able to adapt to a purely dissipative regime, which is fundamentally different from the vortex-shedding scenario. In particular, the results show that:
\begin{itemize}
    \item purely \textbf{reference-guided} training may suffer from limited extrapolation capability when the dynamics evolve beyond the training window;
    \item randomization in the training time window can enhance robustness but may introduce mild small-scale noise;
    \item \textbf{reference-free} reward formulations significantly improve stability while maintaining spectral and energetic accuracy.
\end{itemize}

\subsection{Computational details}
\label{comp-details}
This section summarizes the main computational aspects related to the RL experiments.
In particular, Table \ref{tab:cpu-times} reports the training times for each RL-EF approach and for both test cases, while Table \ref{tab:hyperparameters-rl} lists the main RL hyperparameters for reproducibility.
Regarding the computational times, all the experiments have been performed on the Galileo cluster G100 at CINECA (Italy), on Intel Xeon Platinum 8260 nodes. No GPU acceleration was used. 
Table \ref{tab:cpu-times} highlights that the CPU time needed for the DNS in Test Case 2 is much larger than the CPU time required for the DNS in Test Case 1. This happens due to the larger number of Finite Elements nodes, as can be seen comparing Tables \ref{tab:dofs-cyl} and \ref{tab:dofs-decaying}.
On the one hand, in both the test cases, DD-RL-EF-rand requires less CPU time for each episode in average, with respect to the other approaches. As outlined in Table \ref{tab:acronyms-all}, the DD-RL-EF-rand episodes last $N_{\text{train}}^{\text{rand}}=N/10$ time steps, while in all the other methods each episode has $N_{\text{train}}=N/4$ time steps.
On the other hand, DD-RL-EF-rand requires more CPU time to collect the filtered DNS data, as each episode is initialized randomly within the interval $[0, T_{\text{train}}]$, i.e., we require the DNS data in $[0, T_{\text{train}}+T_{\text{train}}^{\text{rand}}]$.
Additionally, we remark that the reference-free approaches, i.e., DF-RL-EF and SP-DF-RL-EF, do not require any DNS data. This allows to save hours of computational time, while maintaining a comparable, or even more accurate, solution with respect to the reference-guided approaches (e.g., see the global spectrum errors in Figures \ref{fig:glob-errs-cyl} and \ref{fig:glob-errs-decaying}). 
The number of episodes specified in Table \ref{tab:cpu-times} corresponds to the point at which the training reward stabilized, after which the model was saved to perform time extrapolation tests.

\begin{table*}[h]
\centering
\caption{Computational training times for RL-EF approaches.}
\label{tab:cpu-times}
{
\begin{tabular}{ccccc}
\toprule
\textbf{Test case}&{\textbf{Algorithm}} & \textbf{Average CPU time per episode}& \textbf{Episodes} &\textbf{DNS CPU time}
\\ \midrule
\multirow{4}{*}{1} & DD-RL-EF & $15$ min $38$ s & $150$ &$63$ min $41$ s  \\  
&DD-RL-EF-rand & $6$ min $18$ s & $300$ & $89$ min $2$ s \\  
&DF-RL-EF & $15$ min $18$ s & $200$& Not needed \\
\midrule
\multirow{5}{*}{2} & DD-RL-EF & $32$ min $14$ s & $90$& $5$ h $10$ min $30$ s  \\  
&DD-RL-EF-rand & $13$ min $11$ s &$210$ & $7$ h $13$ min $45$ s \\  
&DF-RL-EF & $32$ min $9$ s & $90$ & Not needed \\ 
&SP-DF-RL-EF & $32$ min $14$ s &$60$ & Not needed \\
\bottomrule
\end{tabular}}
\end{table*}

\begin{table*}[h]
\caption{Hyperparameters used for RL training.}
\label{tab:hyperparameters-rl}
\centering
{
{
\begin{tabular}{cc}
\toprule
\textbf{Hyperparameter} & \textbf{Value}\\
\midrule
Discount factor $\gamma$ & {$0.99$}\\
Learning rate & {$\num{1e-5}$}\\
Network architecture & {$2$ hidden layers of $64$ neurons each}\\
Network activation & {{ReLU}}\\
Network optimizer & {Adam} \\
Max gradient norm & {$5$} \\
Batch size & {$128$} \\
Target update interval &  $5 N_{\text{train}}$ \\
\bottomrule
\end{tabular}}}
\end{table*}

\newpage
\section{Conclusions}
\label{sec:conclusions}

In this work, we introduce \textbf{RL-EF}: a reinforcement learning framework for the dynamic selection of the filter radius in the EF stabilization of convection-dominated incompressible flows. The key idea is to replace heuristic, constant-in-time prescriptions of the filter width with a learned, state-dependent control policy that adapts to the evolving flow dynamics.
The proposed RL-EF methodology reformulates parameter tuning as a sequential decision-making problem. The RL agent learns a global control strategy over a prescribed training window and then applies it in time extrapolation. Both reference data-driven (DD-RL-EF) and fully reference data-free (DF-RL-EF) reward formulations were investigated, including a structure-preserving variant incorporating weak energy and enstrophy constraints (if needed by the specific test case).
The numerical assessment on two benchmark problems with fundamentally different dynamics, flow past a cylinder at $Re = \num{1000}$ and decaying homogeneous turbulence at $Re = \num{40000}$, demonstrates several key findings:
\begin{itemize}
    \item \textbf{Stability without over-diffusion.}\\
    The learned policies prevent the numerical blow-up observed in under-resolved noEF simulations, while avoiding the excessive dissipation typical of EF approaches based on a fixed Kolmogorov length scale. The RL-selected filter radius varies non-trivially in time, confirming that no single constant value is optimal across the entire simulation.
    \item \textbf{Spectral consistency across scales.}\\
    RL-EF strategies preserve large-scale structures while correctly controlling the dissipative range. The resulting energy spectra remain close to the reference filtered DNS, including in the inertial range, demonstrating that adaptive filtering can recover the correct balance between energy-containing and small-scale modes.
    \item \textbf{Generalization beyond the training window.}\\
    In the data-driven setting, the agent is trained only on an early-time window where the flow is not yet fully developed. Despite this restricted information, the learned policy generalizes successfully to later times, including regimes characterized by vortex shedding or long-time decay. This highlights the capability of RL to infer transferable control laws rather than merely interpolating reference data.
    \item \textbf{Effectiveness of physics-informed rewards.}\\
    Most importantly, the data-free formulation achieves performance comparable to the data-driven approach. By relying solely on residual control and gradient-based indicators, and not on computationally costly DNS reference, the DF-RL-EF agent autonomously discovers filtering strategies that are both stable and spectrally accurate. This shows that physically meaningful reward design is sufficient to guide the learning process toward robust stabilization policies, even in the absence of reference solutions.
\end{itemize}
Across both test cases, the learned filtering policies adapt to the qualitative nature of the flow regime: oscillatory and highly dynamic in the cylinder wake, more conservative and persistent in decaying turbulence. This behaviour confirms that the RL framework captures problem-dependent stabilization requirements rather than enforcing a uniform dissipative mechanism.
The proposed approach provides a robust and flexible alternative to manual parameter calibration in EF-type methods. It removes the need for fixed tuning rules based on mesh size or Kolmogorov scaling and replaces them with adaptive, state-aware control strategies learned directly from the flow evolution.

Several research directions naturally emerge from the present study. These aspects have not been explored in the present work, as our primary goal was to provide a preliminary study of the robustness of the proposed approach. A more extensive investigation of these directions will therefore be the subject of future studies.
First, the extension to three-dimensional turbulent flows represents a fundamental step toward realistic high-Reynolds-number applications, where the increased scale separation and richer nonlinear interactions pose substantially greater stabilization challenges. 
Second, a particularly promising direction is the development of a space-dependent filtering strategy. In the current formulation, the filter radius is uniform over the computational domain at each time step. However, the optimal amount of stabilization is generally highly localized, especially in flows with shear layers, vortex shedding, or boundary layers. Allowing the filter radius to vary spatially could significantly enhance accuracy. This could be achieved, for instance, through multi-agent reinforcement learning, where different agents control local regions of the domain, or by leveraging graph-based neural architectures capable of operating directly on mesh connectivity and enabling localized control.
Third, the present framework remains grid-dependent, as the observation space coincides with the discrete velocity degrees of freedom. Designing a mesh-independent representation of the state, making use of graphs, for example, would enhance transferability across meshes and resolutions, paving the way toward discretization-agnostic stabilization policies.
Finally, the adoption of a continuous action space would remove the need for discretizing the filter radius and allow finer control of the stabilization mechanism. While preliminary experiments with continuous-action algorithms did not yet yield competitive performance, further investigation of policy-gradient methods and tailored network architectures may unlock additional gains in adaptability and robustness.

The results presented here suggest that reinforcement learning offers a promising pathway toward self-adaptive, physically consistent stabilization techniques for high-Reynolds-number flow simulations.
\section*{Acknowledgments}
AI, MS, and GR thank INdAM-GNCS: Istituto Nazionale di Alta Matematica –– Gruppo Nazionale di Calcolo Scientifico for the support.
MS and GR thank the ``20227K44ME - Full and Reduced order modelling of coupled systems: focus on non-matching methods and automatic learning (FaReX)" project – funded by European Union – Next Generation EU  within the PRIN 2022 program (D.D. 104 - 02/02/2022 Ministero dell’Università e della Ricerca).
GR acknowledges the consortium iNEST (Interconnected North-East Innovation Ecosystem), Piano Nazionale di Ripresa e Resilienza (PNRR) – Missione 4 Componente 2, Investimento 1.5 – D.D. 1058 23/06/2022, ECS00000043, supported by the European Union's NextGenerationEU program.
MS thanks the ECCOMAS EYIC Grant ``CRAFT: Control and Reg-reduction in Applications
for Flow Turbulence".
\bibliographystyle{abbrv}
\bibliography{main}

\end{document}